 \documentclass[final ]{IEEEtran}

 \usepackage{datetime}  
\usepackage{nicematrix} 
\usepackage{multicol}
\usepackage{mathtools}

\IEEEoverridecommandlockouts              

\usepackage{amsmath,amssymb}
\usepackage{amsthm,thmtools}
\usepackage{amsfonts} 
\usepackage{blkarray}

 \usepackage{graphicx}
\usepackage{enumerate}
\usepackage{color}      

\usepackage{verbatim}

\usepackage{amssymb}
\usepackage{amsbsy}
\usepackage{amsmath,lipsum}
\usepackage{setspace}
\usepackage{url}

\usepackage{float}
 \usepackage{dsfont}
 \usepackage{cuted} \stripsep-3pt

 \usepackage{caption} 
\usepackage{todonotes}
\newcommand{\rotodo}[1]{\todo[inline]{\textbf{From RO:} #1}}

 \newcommand{\real}{\operatorname{Re}}
 
 \newcommand{\diag}{\operatorname{diag}}
 \newcommand{\sign}{\operatorname{sgn}}
 \newcommand{\tr}{\operatorname{tr}}
 \newcommand{\adju}{\operatorname{adj}}
 \newcommand{\suchthat}{\, | \,}

\newcommand{\spanop}{\operatorname{span}}

\newcommand*\diff{\mathop{}\!\mathrm{d}}
\newcommand{\trace}{\operatorname{tr}}
\newcommand\scalemath[2]{\scalebox{#1}{\mbox{\ensuremath{\displaystyle #2}}}}

\declaretheorem[name={Example},qed={\lower-0.3ex\hbox{$\square$}} ] {Example}

\declaretheorem[name={Definition}  ] {Definition}
\newtheorem{Theorem}{Theorem}
\newtheorem{Proposition}[Theorem]{Proposition}
\newtheorem{Lemma}[Theorem]{Lemma}
\newtheorem{Corollary}[Theorem]{Corollary}
\declaretheorem[name={Remark}  ] {Remark}

\newcommand {\R}{\mathbb R}

\newcommand {\C}{\mathbb C}

\newcommand{\be}{\begin{equation}}
\newcommand{\ee}{\end{equation}}

\usepackage{xcolor}
\usepackage{soul}

\usepackage{lineno}

\DeclareMathOperator{\argmax}{argmax}

\title{\LARGE \bf
Verifying   $k$-Contraction without
Computing $k$-Compounds \\  \vspace*{20pt} \normalsize  \today{}  }
 
\begin{document}
\author{Omri Dalin, Ron Ofir, Eyal Bar Shalom,  Alexander Ovseevich,
 Francesco Bullo, and Michael~Margaliot 
\thanks{OD, RO, EBS, AO and MM are with the School of Elec.  Eng., Tel Aviv University, Israel 69978. FB is with the Mechanical Engineering Department at UC Santa Barbara. This research was partly supported by a research  grant from  the ISF. The work of FB was supported in part by AFOSR grant FA9550-22-1-0059. Corresponding author: Michael Margaliot (michaelm@tauex.tau.ac.il)}}

%
\maketitle
%

\begin{abstract}
%
Compound matrices have found applications in many fields of science including  systems and control theory. In particular, a sufficient condition for $k$-contraction is  that a logarithmic  norm (also called matrix measure) of the $k$-additive compound of the Jacobian is uniformly negative. However, this may be difficult  to check in practice because the~$k$-additive compound of an~$n\times n$ matrix has dimensions~$\binom{n}{k}\times \binom{n}{k}$. 
For an~$n\times n$ matrix~$A$, we prove a duality relation between the~$k$
 and~$(n-k)$ compounds of~$A$. We use this duality relation to derive a sufficient condition for  $k$-contraction that does not require the  computation of any $k$-compounds. 
 We demonstrate our results by deriving a sufficient condition for~$k$-contraction of an $n$-dimensional Hopfield network
 that does not require to compute any compounds.
In particular, for~$k=2$ this  sufficient condition implies that the network is $2$-contracting and this   
   implies a strong asymptotic property: 
 every bounded  solution of the network  converges to an equilibrium point, that may not be unique. This is   relevant, for example, 
 when using the Hopfield network as an associative memory that stores patterns as equilibrium points of the dynamics. 
%
\end{abstract}
%
%


%
\section{Introduction}

A nonlinear dynamical system  is called contracting if
any two solutions approach one another at an exponential rate~\cite{bullo_contraction}.
This implies many useful asymptotic  properties that resemble those of   asymptotically stable linear systems. For example, 
if the vector field is time-varying and~$T$-periodic and the state-space is convex and bounded 
then the system admits a unique $T$-periodic solution that is globally exponentially  stable~\cite{LOHMILLER1998683,entrain2011,sontag_cotraction_tutorial}. 
If     the   periodicity of the vector field   represents  a $T$-periodic excitation 
then this implies that the system \emph{entrains} to the excitation. 
In fact,  contracting
systems have a well-defined frequency response, as shown in~\cite{freq_convergent} in the closely-related 
context of convergent  systems~\cite{pavlov_book}.

These properties are
important in many applications
ranging from   the synchronized 
response of   biological processes 
to periodic excitations 
like the cell cycle~\cite{RFM_entrain,entrain2011} or 
the 24h solar day, to the entrainment of synchronous generators to the 
frequency of the electric grid.

In particular, if the vector field of a contracting system  is time-invariant then the system admits a   globally exponentially stable equilibrium point. 
Contractivity implies many other useful properties e.g., a contractive system is input-to-state stable~\cite{desoer_1972,sontag2010contractive}.

An important advantage of
contraction theory is that  there exists  a simple  sufficient condition for contraction, namely, 
that  a logarithmic norm (also called matrix  measure) of the Jacobian
of the vector field is uniformly negative.
 For the~$L_1$, $L_2$, and~$L_\infty$ norms this
 sufficient condition is easy to check,
 and, in particular,  does not require explicit knowledge of the solutions of the system. 

Contraction theory has found numerous applications in
robotics~\cite{robo_contraction},
synchronization in
multi-agent systems~\cite{slotine_sync}, 
the design of observes and closed-loop controllers~\cite{control_syn_vincent,andrieu_trans}, neural networks and learning theory~\cite{slotine_PLOS_ONE},
and more.

    However, many systems cannot be studied using contraction theory. For example, if the dynamics admits more than a single equilibrium then the system is clearly not contracting. Existence of more than a single equilibrium, i.e., multi-stability, is prevalent 
    in many important  mathematical models and real-world systems.  Ecological models that include several equilibrium points allow switching 
    between several possible behaviours e.g., outbreaks~\cite{bistab_in_eco}.
    Epidemiological models typically include at least two equilibrium points corresponding to the disease-free steady state and the endemic steady state~\cite{kiss_epi_2017}. Multi-stability 
    in biochemical and cellular systems allows to transform graded signals
    into an all-or-nothing response 
    and to ``remember''  transitory stimuli~\cite{multi_stable_1999,multi_stab_sontag}. 
Other important examples of systems that are not contractive are systems whose trajectories contract at a rate that is slower than exponential, systems that are almost globally stable,
and more~(see, e.g.,~\cite{Angeli2021nonosc}).

There is considerable  interest in extending contraction theory to 
systems that are not contractive in the usual sense, see, e.g.~\cite{weak_contraction,weak_contra_bullo}.
Motivated by the seminal work of Muldowney~\cite{muldowney1990compound}, Wu et  al.~\cite{kordercont} recently introduced the notion of $k$-contractive systems. 
Roughly speaking, the solutions of these systems contract~$k$-dimensional parallelotopes.
For~$k=1$, this reduces to standard contraction. However,~$k$-contraction with~$k>1$ can  be used to analyze the asymptotic behaviour of systems that are not contractive. For example, every bounded solution of a time-invariant   $2$-contractive system converges 
to an equilibrium, that is not necessarily unique~\cite{li1995}. 

A sufficient condition for  $k$-contraction is that a log norm  of  the~$k$-additive compound of the Jacobian of the vector field is uniformly negative.
However,  this is not always easy to check, 
as the~$k$ compound of an~$n\times n$ matrix has dimensions~$r\times r$, with~$r:=\binom{n}{k}$. For more on the applications of compound systems to
systems and control theory, see e.g.~\cite{wu2021diagonal, margaliot2019revisiting, ofir2021sufficient, ron_DAE,grussler2022variation,LI1999191} and the recent tutorial~\cite{comp_long_tutorial}.

Here,  we derive  new duality relations between the $k$ and the~$(n-k)$ compounds  of a matrix~$A \in \C^{n\times n}$.
  We then use these duality relations to 
  derive a new  sufficient condition for
$k$-contraction that does \emph{not}
require computing any compounds. We demonstrate our theoretical results by deriving a sufficient condition for $k$-contraction in an~$n$-dimensional Hopfield neural network.  This system typically admits  more than a single equilibrium   and is thus not contractive (i.e., not~$1$-contractive) with respect to~(w.r.t.) any norm. Our condition does not require to compute any compounds. For~$k=2$, this condition  implies a strong asymptotic property:   any bounded solution of the network  converges to an equilibrium.

The remainder of this paper   is organized as follows. The next section reviews known definitions and results that are used later on. Section~\ref{sec:main} derives   the  duality relations between compounds.
Section~\ref{sec:comp_free} shows how
these duality relations can be used to prove~$k$-contraction without computing any compounds. 
We demonstrate the usefulness of the theoretical results by deriving a sufficient condition for~$k$-contraction in Hopfield neural networks. 
The final section concludes. 

We use standard notation. 
Vectors [matrices] are denoted by small [capital] letters. For a matrix~$A$, $a_{ij}$ is entry~$(i,j)$ of~$A$, and~$A^T$ is the transpose of~$A$. For a square matrix~$A$, $\tr(A)$ [$\det(A)$] is the trace [determinant]  of~$A$. 
For a symmetric matrix~$Q\in\R^{n\times n}$, we use~$Q\succ 0$ [$Q\succeq 0$] 
to denote that~$Q$ is positive-definite [non-negative definite]. A square matrix is called
  \emph{anti-diagonal}   if all its
  entries are zero except those on the diagonal going from the lower left corner to the upper right corner, known as the \emph{anti-diagonal}. 

\section{Preliminaries}
This section reviews  several known definitions  and results that are used later on. 
Let~$Q(k,n)$ denote all
 the~$\binom{n}{k}$
increasing sequences of~$k$ integers from the set~$\{1,\dots,n\}$, ordered lexicographically.
For example,
\[
Q(3,4)= ( (1,2,3), (1,2,4), (1,3,4), (2,3,4) ).
\]
Let~$A\in\C^{n\times m }$.   Fix~$k\in\{1,\dots,\min(n,m) \}$.
For two sequences~$\alpha  \in Q(k,n)$,  $\beta \in Q(k,m)$, let~$A[\alpha|\beta]$ denote the~$ k\times k$ submatrix   obtained by taking the entries of~$A$ in  the rows indexed by~$\alpha$ and the columns indexed by~$\beta$. For example
\[
A[ (2,4)|(1,2)]=\begin{bmatrix}
a_{21} & a_{22}\\
a_{41} &a_{42}
\end{bmatrix} .
\]
The  minor of~$A $ corresponding to~$\alpha ,  \beta$  is
\[
A(\alpha|\beta):=\det (A[\alpha|\beta])  .
\]
For example, if~$m=n$ then~$Q(n,n)$ includes the single element~$\alpha=(1,\dots,n)$, so~$A[\alpha|\alpha]=A$,
and~$A(\alpha|\alpha)=\det(A)$.

 \subsection{Compound Matrices}
The~$k$-multiplicative compound of a matrix~$A$ is a   matrix that collects  all the~$k$-minors of~$A$. 
\begin{Definition}\label{def:multi} 
%
Let $A\in \C^{n\times m}$ 
and fix $k \in \{ 1,\dots,\min  (n,m) \}$. 
The \emph{$k$-multiplicative compound} of~$A$, denoted~$A^{(k)}$, is the~$\binom{n}{k}\times \binom{m}{k}$
matrix that contains all the~$k$-minors of~$A$  ordered lexicographically. 
\end{Definition}

For example, if~$n=m=3$ and~$k=2$ then
\[
A^{(2)}= \begin{bmatrix}
A((12)|(12)) & A((12)|(13)) & A((12)|(23)) \\ 
A((13)|(12)) & A((13)|(13)) & A((13)|(23)) \\ 
A((23)|(12)) & A((23)|(13)) & A((23)|(23)) 
\end{bmatrix}.
\]
In particular, Definition~\ref{def:multi} 
implies that~$A^{(1)}=A$,
and if~$n=m$ then~$A^{(n)}=\det(A)$.
Note also  that by definition~$(A^T)^{(k)}=(A^{(k)})^T$. In particular,  if~$A$ is symmetric 
then~$(A^{(k)})^T =(A^T)^{(k)}= A^{(k)} $, so~$A^{(k)}$ is also symmetric.

The term multiplicative compound is justified  by the following important result. 
\begin{Theorem}[Cauchy-Binet Theorem]\label{thm:CB}
 Let~$A\in\C^{n \times m}$ and~$B \in \C^{m\times p}$. 
 Then for any~$k\in\{1,2,\dots,\min\{ n,m,p \}\}$, we have
 \[
 (AB)^{(k)}= A^{(k)}B^{(k)}.
 \]
\end{Theorem}

Note that for~$m=p=k=n$ this reduces to~$\det(AB)=\det(A)\det(B)$. 

Thm.~\ref{thm:CB} implies in particular that
if~$A$ is square and invertible then~$A^{(k)} $ is also invertible, and~$(A^{(k)})^{-1} =(A^{-1})^{(k)} $. 

If~$A=\diag(\lambda_1,\lambda_2,\lambda_3)$ then Definition~\ref{def:multi} implies that~$A^{(2)}=\diag(\lambda_1\lambda_2, \lambda_1 \lambda_3,\lambda_2\lambda_3)$, so
every eigenvalue of~$A^{(2)}$ is a product of two eigenvalues of~$A$.
More generally, 
the multiplicative compound has a useful spectral property. Let~$\lambda_1,\dots,\lambda_n$ denote the eigenvalues of~$A\in\C^{m\times m}$. Then the~$\binom{n}{k}$ eigenvalues of~$A^{(k)}$ are all the  products of $k$ eigenvalues of~$A$, i.e.
\begin{equation}\label{eq:spec_multi}
    \prod_{i=1}^k \lambda_{\alpha_i},\text{for all } \alpha \in Q(k,n).
\end{equation}
 A similar property also holds for the singular values of~$A^{(k)}$.
Eq.~\eqref{eq:spec_multi} implies the Sylvester-Franke identity:
\be\label{eq:sylv}
\det(A^{(k)})=(\det(A))^{\binom{n-1}{k-1}} 
\ee
(see, e.g.,~\cite{Sylvester-Franke}).

The $k$-multiplicative compound can be used to 
study  the evolution of 
the volume of~$k$-dimensional parallelotopes under a differential equation.
The parallelotope
with vertices~$x^1,\dots,x^k \in \R^n$
(and the zero vertex)  is 
\[
P(x^1,\dots,x^k): =  \{\sum_{i=1}^k c_i x^i \suchthat  c_i\in[0,1]\}.
\]
Let~$X:=\begin{bmatrix} x^1 & \dots &  x^k \end{bmatrix}^{(k)}$. 
Note that~$X$ has dimensions~$\binom{n}{k}\times 1$.
It is well-known (see, e.g.,~\cite{comp_long_tutorial})  that the volume of~$P (x^1,\dots,x^k)$ satisfies 
\[
\text{volume}(P(x^1,\dots,x^n)) = |\begin{bmatrix} x^1 & \cdots &  x^k \end{bmatrix}^{(k)} |_2 ,
\]
where~$|\cdot|_p$ denotes  the $L_p$ norm. 
For the particular case~$k=n$ this becomes the well-known formula
\[
\text{volume}(P(x^1,\dots,x^n))=|\det(x^1,\dots,x^n)|.
\]
  To study the time evolution of such $k$-volumes under the action of   a differential equation 
  requires the~$k$-additive compound of a square matrix.
\begin{Definition}\label{def:add_comp}
%
Let~$A\in\C^{n\times n}$.
The \emph{$k$-additive compound}
matrix of~$A$ is the~$\binom{n}{k}\times \binom{n}{k}$ matrix defined by:
\begin{align}\label{eq:A^[k]:=ddeps}
    A^{[k]} := \frac{d}{d\varepsilon} (I_n+\varepsilon A)^{(k)} |_{\varepsilon=0} . 
\end{align}
%
\end{Definition}
The derivative here is well-defined, as every entry of~$(I_n+\varepsilon A)^{(k)}$ is a polynomial in~$\varepsilon$.
Note that this definition implies that
\be\label{eq:akasexp}
A^{[k]} =  \frac{d}{d\varepsilon} (\exp(A\varepsilon ))^{(k)} |_{\varepsilon=0} ,
\ee
and also the Taylor series expansion 
\begin{align}\label{eq:(I+epsilonA)^k}
    (I_n+\varepsilon A)^{(k)} = I_r + \varepsilon A^{[k]} + o(\varepsilon),
\end{align}
%
where~$r:=\binom{n}{k}$.

 Definition~\ref{def:add_comp} implies that many properties of~$A^{[k]}$ can be deduced from properties of~$A^{(k)}$. For example,~\eqref{eq:spec_multi} implies that
 the~$\binom{n}{k}$ eigenvalues of~$A^{[k]}$ are
 \begin{equation}\label{eq:spec_add}
    \sum_{i=1}^k \lambda_{\alpha_i},\text{for all } \alpha \in Q(k,n).
\end{equation}
Also,~\eqref{eq:sylv} implies that 
 \be\label{eq:traceak}
 \tr(A^{[k]})=\binom{n-1}{k-1}\tr(A). 
 \ee

The next  result provides 
a useful explicit formula for~$A^{[k]}$ in terms of the entries~$a_{ij}$ of~$A$.
  Recall that any entry of~$A^{(k)}$ is a minor~$A(\alpha|\beta)$.
Thus, it is natural to index the entries of~$A^{(k)}$  and~$A^{[k]}$  using~$\alpha,\beta \in Q(k,n)$. 
%
\begin{Proposition}\label{prop:Explicit_A_k}\cite{schwarz1970}
%
Fix~$\alpha,\beta \in Q(k,n)$ and let~$\alpha=\{i_1,\dots,i_k\}$ and~$\beta=\{j_1,\dots,j_k\}$. Then the entry of~$A^{[k]}$ corresponding to~$(\alpha,\beta)$ is
equal to:
\begin{enumerate}
    \item $\sum_{\ell=1}^{k} a_{ i_{\ell} i_{\ell} }$, if  $i_{\ell} = j_{\ell}$  for all $\ell \in \{ 1,\hdots,k \}$;
    \item $(-1)^{\ell +m}
    a_{i_{\ell} j_{m}} $, 
      if all the indices  in $ \alpha  $  and $ \beta$   agree,  except  for   
   a single index $ i_{\ell} \ne j_m$; and
    \item $0$, otherwise.
\end{enumerate}
%
\end{Proposition}
 Note that the   first case in the proposition corresponds to the diagonal entries of~$A^{[k]}$. Also,  the proposition
implies in particular that~$A^{[1]}=A$, and~$A^{[n]}= \sum_{\ell=1}^{n} a_{\ell \ell}= \tr(A)$.

\begin{Example}\label{ex:A^[3]}
%
For~$A\in\R^{3\times 3}$ and~$k=2$,
Prop.~\ref{prop:Explicit_A_k} yields
\[ A^{[2]} =
\begin{pNiceMatrix}[first-row,  last-col]
     (1,2) & (1,3) & (2,3)        \\
  a_{11}+a_{22}   & a_{23}  & -a_{13}     & (1,2) \\
  a_{32}      & a_{11}+a_{33} &a_{12}     & (1,3) \\
  -a_{31}  & a_{21}  & a_{22}+a_{33}      & (2,3)
\end{pNiceMatrix},
\]
where the indexes~$\alpha \in Q(2,3)$ [$ \beta \in Q(2,3)$] are marked on  right-hand [top] side  
  of  the matrix.
For example, the entry in the second row and 
third column
of $A^{[3]}$ corresponds to 
$(\alpha | \beta) = ((1,3 ) |  ( 2,3  ) )$.
As~$\alpha$ and~$\beta$ agree in all indices except
for the  index
$i_{ 1} =1$ and $j_{ 1} =2$,
this  entry is equal to~$  (-1)^{1+1} a_{12} =  a_{12}$.
%
\end{Example}


Prop.~\ref{prop:Explicit_A_k} implies that the mapping~$A\to A^{[k]}$ is linear and, in particular, 
\[
(A+B)^{[k]}=
A^{[k]}+B^{[k]},
\text{ for any } A,B\in \C^{n\times n}. \]
This justifies the term additive compound.

It is useful  to consider the additive compound of a matrix under a coordinate transformation.  Let~$A \in \C^{n \times n}$. Fix~$k \in \{1,\dots,n\}$, and 
an invertible matrix~$T \in \C^{n \times n}$. Then
\begin{equation}\label{eq:add_T}
    (TAT^{-1})^{[k]} = T^{(k)}A^{[k]}(T^{(k)})^{-1}.
\end{equation}

To explain the use of the additive compound to study $k$-contraction in nonlinear dynamical systems, we briefly review   log  norms (also called matrix measures and Lozinski\v{i} meaures). 
\subsection{Logarithmic norms}
A   norm $|\cdot |:\R^{n}\to\R_{+}$ induces a matrix norm $\|\cdot \|:\R^{n\times n}\to\R_{+}$ defined 
by~$\|A\|:=\max_{|x|=1} |Ax|$,
and a log norm~$\mu  :\R^{n\times n}\to\R$   defined by 
    \[
    \mu (A) :=\lim_{\varepsilon \downarrow 0} \frac{\|I+\varepsilon A\|-1}{\varepsilon}.
    \]
It is well-known that  the solution of~$\dot x=Ax$ satisfies
\[
\frac{d}{dt} \log(|x(t)|)\leq \mu (A),
\]
where the derivative here is the upper right Dini derivative.
    
Log  norms play an important role in numerical linear algebra and in contraction theory   (see e.g.,~\cite{sontag_cotraction_tutorial, vidyasagar2002nonlinear,strom1975logarithmic}). For the~$L_1$, $L_2$, and~$L_\infty$ norms, there exist closed-form expressions for the induced log norms, namely, 
    \begin{align*}
        \mu_1(A) &= \displaystyle \max_{j}\Big(a_{jj}+\sum^n_{\substack{i=1\\
                  i\neq j}}|a_{ij}|\Big),               \\
        \mu_2(A) &= \lambda_1\Big(\frac{A+A^T}{2}\Big), \\
        \mu_{\infty}(A) &= \displaystyle \max_{i}\Big(a_{ii}+\sum^n_{\substack{j=1\\
                  i\neq j}}|a_{ij}|\Big),
    \end{align*}
where~$\lambda_i(S)$ denotes the~$ith$ largest eigenvalue of the symmetric matrix~$S$, that is,~$\lambda_1(S)\geq\lambda_2(S)\geq\hdots\geq\lambda_n(S)$.  Using Prop.~\ref{prop:Explicit_A_k},
this can be generalized to 
closed-form expressions for the   induced log norms of the additive compounds of a matrix.
\begin{Proposition}(see  e.g.,~\cite{muldowney1990compound})\label{prop:k_add_meas}
Let~$A\in\R^{n\times n}$, and fix~$k\in\{1,\dots,n\}$. Then
    \begin{align*}
        \mu_1(A^{[k]}) &= \max_{\alpha \in Q(k,n)} \sum_{i \in \alpha} \left(a_{ii} + \sum_{j \notin \alpha} |a_{ji}|\right), \\
        \mu_2(A^{[k]}) &= \sum_{i=1}^k\lambda_i\Big(\frac{A+A^T}{2}\Big), \\
        \mu_\infty(A^{[k]}) &= \max_{\alpha \in Q(k,n)} \sum_{i \in \alpha} \left(a_{ii} + \sum_{j \notin \alpha} |a_{ij}|\right).
    \end{align*}
\end{Proposition}
 
It is straightforward to verify that if~$H\in\R^{n\times n}$ is invertible  then the scaled  norm~$|\cdot|_H$, defined by~$|x|_H:=|Hx|$, induces the log norm 
\begin{equation}\label{eq:weight_mat_meas}
    \mu_H(A):=\mu(H A H^{-1}).
\end{equation}
  $L_p$ norms are invariant under permutations and sign changes, i.e.,
  if~$P$ [$S$] is a 
  a permutation [signature] 
  matrix then~$|x|_p=|PSx|_p$ for any~$x$. This yields  the following result.
\begin{Lemma}\label{lem:pnorm_U_inv}
  Let $\mu_p(\cdot)$ denote the   log norm  induced by~$|\cdot|_p$. If~$U \in \R^{n \times n}$ is the product of a permutation matrix and a signature matrix, then
\[
        \mu_{p,U}(A) = \mu_p(A),  \text{ for any } A\in\R^{n\times n}.
\]
\end{Lemma}

We also require  a duality result for the log norm 
that follows from  a well-known relation between an~$L_p$ norm and its dual norm.
\begin{Lemma}\label{lem:mat_meas_dual}
    Let $p,q \in [1,\infty]$ such that $ {p}^{-1} +  {q} ^{-1}= 1$. Then
\[
        \mu_p(A) = \mu_q(A^T), 
   \text { for any }
   A \in \R^{n \times n}.
    \]
\end{Lemma}
\begin{IEEEproof}
   We begin by proving a duality for the induced matrix norm. Using the definition of the induced matrix norm and the dual norm of $L_p$ norms, we have
    \begin{align*}
        \|A\|_p &= \max_{|x|_p = 1} |Ax|_p \\
        &= \max_{|x|_p = 1} \max_{|y|_q = 1} |(Ax)^Ty| \\
        &= \max_{|x|_p = 1} \max_{|y|_q = 1} |x^T (A^Ty)| \\
        &\le \max_{|x|_p = 1} \max_{|y|_q = 1} |x|_p |A^Ty|_q \\
        &= \|A^T\|_q,
    \end{align*}
    where we used H{\"o}lder's inequality. 
    Since $(A^T)^T = A$, this implies that $\|A\|_p = \|A^T\|_q$.
   Thus, 
    \begin{align*}
        \mu_p(A) &= \lim_{h \to 0^+} \frac{\|I_n + hA\|_p - 1}{h} \\
        &= \lim_{h \to 0^+} \frac{\|I_n + hA^T\|_q - 1}{h} \\
        &= \mu_q(A^T),
    \end{align*}
    and this completes the proof.
\end{IEEEproof}
\subsection{ $k$-contraction  }

Motivated by the seminal work of Muldowney~\cite{muldowney1990compound}, Wu et  al.~\cite{kordercont}
recently introduced the notion of $k$-contractive systems. Roughly speaking, the solutions of these systems contract~$k$-dimensional parallelotopes.
For~$k=1$, this reduces to standard contraction. 

Consider the   time-varying non-linear  system:
\begin{equation}\label{eq:sys}
    \dot x = f(t,x),
\end{equation}
where $x \in \Omega \subseteq \R^n$,
and~$\Omega$ is a convex set.
We assume that~$f $ is~$C^1$,  and denote its Jacobian  with respect to~$x$ by~$J(t,x) := \frac{\partial  }{\partial x}f(t,x)$.
A sufficient condition for~$k$-contraction with rate~$\eta>0$ is that
\be\label{eq:infik}
\mu(J^{[k]}(t,x))\leq-\eta<0, \text{ for all } x\in \Omega, t\geq 0 . 
\ee
For~$k=1$, this reduces to the
standard
sufficient condition for contraction. However, for~$k>1$ this condition is weaker than the one required for~$1$-contraction.
As a simple example, 
a matrix~$A \in\R^{n\times n}$ is Hurwitz iff there exists~$P\succ 0$ such that~$PA+A^TP \prec 0$, that is, iff 
$\mu_{2,P^{1/2}}(A)<0$~\cite{sontag_cotraction_tutorial}.
This implies that~$A^{[2]} $ is
contractive  w.r.t. some  scaled~$L_2$ norm iff~$A^{[2]}$
is Hurwitz, that is, iff 
the sum of any two eigenvalues of~$A$ has a negative real part.
Note that this spectral property implies that any bounded solution of~$\dot x=Ax$ converges to the origin. 

$k$-contraction has several important implication. First, every bounded solution of a~$2$-contractive time-invariant nonlinear dynamical system converges to an equilibrium (that may not be unique)~\cite{li1995}.   Second,
for LTV systems  $k$-contraction implies the existence of a stable subspace. 
\begin{Proposition} \cite{muldowney1990compound} 
\label{prop:subspace}  
	Suppose that the LTV  system~$\dot x(t)=A(t)x(t)$,
	where~$A(t)$ is a continuous matrix function,
	is uniformly stable. Let~$x(t_0) =x_0$ denote an initial condition at time~$t_0$. Fix~$k\in\{1,\dots,n\}$.  
The following   two  conditions are equivalent. 
	\begin{enumerate}[(a)]
	\item \label{cond:nk1} The LTV system
	admits an~$(n-k+1)$-dimensional linear subspace 
	$  \mathcal{X} (t_0) \subseteq \R^n$ such that
	\begin{equation} \label{eq:a0x}
	\lim_{t \to \infty }x(t, t_0, x_0) = 0 \text{ for any } x_0 \in \mathcal{X}(t_0) .
	\end{equation}
	\item \label{cond:sysmkop} Every solution of  
	\begin{equation} \label{eq:kinsys}
	\dot{y}(t) = A^{[k]}(t)y(t)  
	\end{equation}
	converges to the origin as~$t\to\infty$.
\end{enumerate}
\end{Proposition}
Note that if~$\mu(A^{[k]}(t))\leq-\eta<0$ then clearly condition~\eqref{cond:sysmkop}  holds, and thus~\eqref{cond:nk1} holds.

\begin{Example}
Consider the LTV~$\dot x(t)=A(t)x(t)$, $x(t_0)=x_0$,  with
\[
A(t)=(1/2)\begin{bmatrix}
-3+3\cos^2(t) & 2-3\cos(t)\sin(t)\\
-2-3\cos(t)\sin(t) & -3+3 \sin^2(t)
\end{bmatrix}.
\]
It can be verified that~$x(t)=\Phi(t,t_0)x_0$, 
where the transition matrix~$\Phi(t,t_0)$ is
\begin{align}\label{eq:ltvtra}
\Phi(t,t_0) & =\begin{bmatrix}
 \cos(t) & \sin(t) \\ -\sin(t)& \cos(t) 
\end{bmatrix}
\diag\left (1,\exp(-3(t-t_0)/2) \right )
\nonumber\\
&\times
\begin{bmatrix}   
 \cos(t_0) & -\sin(t_0) \\ 
 \sin(t_0)& \cos(t_0) 
\end{bmatrix}.
\end{align}
This implies that the LTV is uniformly stable, but not contractive. 
Here,
\[
A^{[2]}(t)=\tr(A(t))=-3/2,
\]
so the LTV is~$2$-contractive, and Prop.~\ref{prop:subspace} 
implies that 
the LTV  
	admits a~$1$-dimensional linear subspace 
	$\mathcal{X} \subseteq \R^2$ such that~\eqref{eq:a0x} holds. 
Indeed, it follows from~\eqref{eq:ltvtra} that~$\spanop(\begin{bmatrix}
\sin(t_0)   \\
\cos(t_0) 
\end{bmatrix})$
is such a subspace. 
\end{Example}

In principle, verifying that~\eqref{eq:infik} holds can be done by first computing $J^{[k]}(t,x)$, for~$p\in\{1,2,\infty\}$, and then 
using 
the expressions in Prop.~\ref{prop:Explicit_A_k}.
However, in practice this is  non-trivial, as~$|Q(k,n)|=\binom{n}{k}$.

\section{Main Results}\label{sec:main}
From here on we fix an integer~$n>0$, and an integer~$k\in\{1,\dots,n-1\}.$
Let~$r:=\binom{n}{k}$.
Note that  the  matrices~$A^{(k)}$, $A^{(n-k)}$, $A^{[k]}$, and~$A^{[n-k]}$ 
all have the same dimensions, namely,~$r\times r$. Our goal is to derive certain duality relations between these matrices, and then use them to relate~$\mu(A^{[k]})$ and~$\mu(A^{[n-k]})$. 

We begin by  defining  an anti-diagonal  matrix~$U(k,n)$  that will be used in the  results below. Denote the lexicographically  ordered sequences in~$Q(k,n)$ by~$\alpha^1,\dots,\alpha^r$.
The signature of~$\alpha^j$ is~$s(\alpha^j):= (-1)^{\alpha^j_1+\dots+\alpha^j_k }$, and the complement of~$\alpha^j$ is
\[
\overline{\alpha^j}:=\{1,\dots,n\}\setminus\alpha^j.
\]
For simplicity, we use set notation here, but we always assume that 
the entries of~$\overline{\alpha^j}$ are arranged in the lexicographic  order. 

\begin{Definition}\label{def:U_matrices}
Let~$U=U(k,n) \in\{-1,0,1\}^{r\times r}$ be the anti-diagonal matrix with entries:
\begin{align}\label{eq:U_r}
   u_{ij}=
    \begin{cases}
    s(\alpha^j), & \text{if } i+j=r+1,\\
    0, & \text{otherwise}.
    \end{cases}
\end{align}
\end{Definition}

\begin{Example}\label{eax:n4k2}
For~$n=4$ and~$k=2$, we have~$\alpha^1= (1,2)$,
$\alpha^2= (1,3)$,
$\alpha^3= (1,4)$,
$\alpha^4= (2,3)$,
$\alpha^5= (2,4)$,
and~$\alpha^6= (3,4)$, so
$s(\alpha^1)=-1$,
$s(\alpha^2)=1$,
$s(\alpha^3)=-1$,
$s(\alpha^4)=-1$,
$s(\alpha^5)=1$,
$s(\alpha^6)=-1$. Thus,  
\begin{align*}
    U&=\begin{bmatrix}
0 & 0 &0 & 0 & 0 &s(\alpha^6) \\
0 & 0 &0 & 0 & s(\alpha^5) &0 \\
0 & 0 &0 & s(\alpha^4) & 0 &0 \\
0 & 0 &s(\alpha^3) & 0 & 0 &0 \\
0 & s(\alpha^2) &0 & 0 & 0 &0 \\
s(\alpha^1) & 0 &0 & 0 & 0 &0  
\end{bmatrix}
\\&=
\begin{bmatrix}
0 & 0 &0 & 0 & 0 &-1 \\
0 & 0 &0 & 0 & 1 &0 \\
0 & 0 &0 & -1 & 0 &0 \\
0 & 0 &-1 & 0 & 0 &0 \\
0 & 1 &0 & 0 & 0 &0 \\
-1 & 0 &0 & 0 & 0 &0  
\end{bmatrix}.
\end{align*}
\end{Example}

\begin{Example}
For~$k=1$, we have~$Q(1,n)=((1),(2),\dots,(n))$, so~$s(\alpha^j)=(-1)^j$, and thus in this particular case the definition of~$U$ yields
\begin{align}\label{eq:simpleu}
    u_{ij}=
    \begin{cases}
    (-1)^j, & \text{if } i+j=r+1,\\
    0, & \text{otherwise}.
    \end{cases}
\end{align}
However, Example~\ref{eax:n4k2} shows that this expression does not hold in general. 
\end{Example}

Note that~$U^TU=U U^T = I_r$.
A direct calculation shows 
that for any~$B \in \C^{r\times r}$ and any~$i,j\in\{1,\dots,r\}$, we have
\begin{align}\label{eq:utu}
(  U^T B U )_{ij}&=s(\alpha^i)  s(\alpha^{r+1-j}) 
b_{r+1-i,r+1-j}\nonumber\\
& 
=s(\alpha^i)  s(\alpha^{j}) 
b_{r+1-i,r+1-j}.
\end{align}

\subsection{Duality between multiplicative compounds}
The next result describes a duality relation between the two multiplicative compound matrices~$A^{(k)}$ and~$A^{(n-k)}$.
\begin{Theorem}\label{thm:dual_,ult}
Fix~$A\in \C^{n\times n}$,
and let~$U\in\{-1,0,1\}^{r\times r}$ be the anti-diagonal
matrix defined in Definition~\ref{def:U_matrices}.
Then
\begin{align}\label{eq::A^(k)_A^(n-k)_alt}
 (A^{(k)})^T 
 U^TA^{(n-k)}U=\det(A)I_r   . 
\end{align}
\end{Theorem}
In other words, $ U^TA^{(n-k)}U$
is the adjugate    matrix of~$(A^{(k)})^T$.
For~$k=1$, $U$ becomes the matrix in~\eqref{eq:simpleu}, and
Eq.~\eqref{eq::A^(k)_A^(n-k)_alt} 
becomes~$A^T 
 U^TA^{(n-1)}U=\det(A)I_n$, which is just~$\adju(A) A =\det(A)I_n$. 
Formulas that are equivalent to~\eqref{eq::A^(k)_A^(n-k)_alt} are known,  see e.g.,~\cite[p. 29]{matrx_ana} (where it appears without a proof),
but without the explicit expression of the matrix~$U$. 

The proof of Thm.~\ref{thm:dual_,ult}
uses  two  auxiliary results. The first result describes a duality relation between~$Q(k,n)$ and~$Q(n-k,n)$.
\begin{Lemma}\label{lemm:qknqqkn}
If~$Q(k,n)=(\alpha^1,\dots,\alpha^r)$ then
\be\label{eq:qnkdi}
Q(n-k,n)=( \overline{\alpha^r},\dots,\overline{\alpha^1} ).
\ee
\end{Lemma}
\begin{IEEEproof}
It is clear that~$Q(n-k,n)$ 
includes the sequences~$\overline{\alpha^i}$, $i\in\{1,\dots,r\}$. The ordering in~\eqref{eq:qnkdi} follows from the fact that the lexicographic ordering of~$Q(k,n)$ 
is~$\alpha^1,\dots,\alpha^r$,
and the definition of the lexicographic ordering.
\end{IEEEproof}

\begin{Example}
 For~$n=4$ and~$k=3$, we have~$r=4$, and 
 \[
    Q(3,4)=((1,2,3),(1,2,4),(1,3,4),(2,3,4)).
\]
Clearly, 
\[
    Q(1,4)=((1),(2),(3),(4)),
\]
and this agrees with~\eqref{eq:qnkdi}.
\end{Example}

 \begin{Lemma}\label{prop:jacobi-compound_equivalence}
For any~$i,j\in\{1,\dots,r\}$, we have 
\be \label{eq:lemma1}
(U^T A^{(n-k)} U)_{ij}=  s(\alpha^i)s(\alpha^j) A(\overline{\alpha^i}|\overline{\alpha^j}).
\ee
\end{Lemma}

\begin{IEEEproof}
By Lemma~\ref{lemm:qknqqkn}, 
\begin{align*}
   ( A ^{(n-k)})_{pq}=A(\overline{\alpha^{r+1-p}}|\overline{\alpha^{r+1-q}}),
\end{align*} for any~$p,q\in\{1,\dots,r\}$.
 Combining this with~\eqref{eq:utu}  
 yields~\eqref{eq:lemma1}. 
\end{IEEEproof}

We can now prove Thm.~\ref{thm:dual_,ult}. We assume that~$A$ is non-singular. The general case follows by a   continuity argument. 
Denote~$Z:= (A^{(k)})^TU^TA^{(n-k)} U $.
Fix~$i,j\in\{1,\dots,r\}$.
Then~$z_{ij}$ is the product of row~$i$ of $(A^{(k)})^T$ and column~$j$ of~$U^TA^{(n-k)} U$, and combining this with Lemma~\ref{prop:jacobi-compound_equivalence} gives 
\begin{align}\label{eq:sij1}
z_{ij} &=  \sum_{\ell=1}^r
A(\alpha^\ell|\alpha^i)
s(\alpha^\ell) s(\alpha^j) A(\overline{\alpha^\ell}|\overline{\alpha^j} ) . 
\end{align}
   Jacobi's Identity~\cite[p.~166]{muir1933}
 asserts that for any~$p,q$:
\begin{equation*} 
    A(\alpha^p|\alpha^q)=\det(A)  s(\alpha^p) 
    s(\alpha^q) 
    B (\overline{\alpha^q}|\overline{\alpha^p}) , 
\end{equation*}
where~$B:=A^{-1}$. In particular,
\begin{equation*} 
    A(\overline{\alpha^\ell}|\overline{ \alpha^j})=\det(A)  s(\alpha^\ell) 
    s(\alpha^j) 
    B ( {\alpha^j}| {\alpha^\ell}),
\end{equation*}
and substituting this in~\eqref{eq:sij1}
yields
\begin{align*}
z_{ij} &= \det(A) \sum_{\ell=1}^r
A(\alpha^\ell|\alpha^i)
s(\alpha^\ell) s(\alpha^j)
s(\alpha^\ell) 
    s(\alpha^j) 
    B ( {\alpha^j}| {\alpha^\ell})\\
    &=\det(A) \sum_{\ell=1}^r
A(\alpha^\ell|\alpha^i)
    B ( {\alpha^j}| {\alpha^\ell}).
\end{align*}
In other words,~$z_{ij}/\det(A)$ is the product of row~$j$ of~$(A^{-1})^{(k)}$ with
column~$i$ of~$A^{(k)}$, that is,
\begin{align*}
    Z^T &=\det(A) (A^{-1})^{(k)} A^{(k)}\\
    &=\det(A) (A^{-1}A)^{(k)} \\
    &=\det(A) I_r ,
\end{align*}
 and this completes the proof of Thm.~\ref{thm:dual_,ult}. 

\begin{figure}
 \begin{center}
  \includegraphics[scale=0.45]{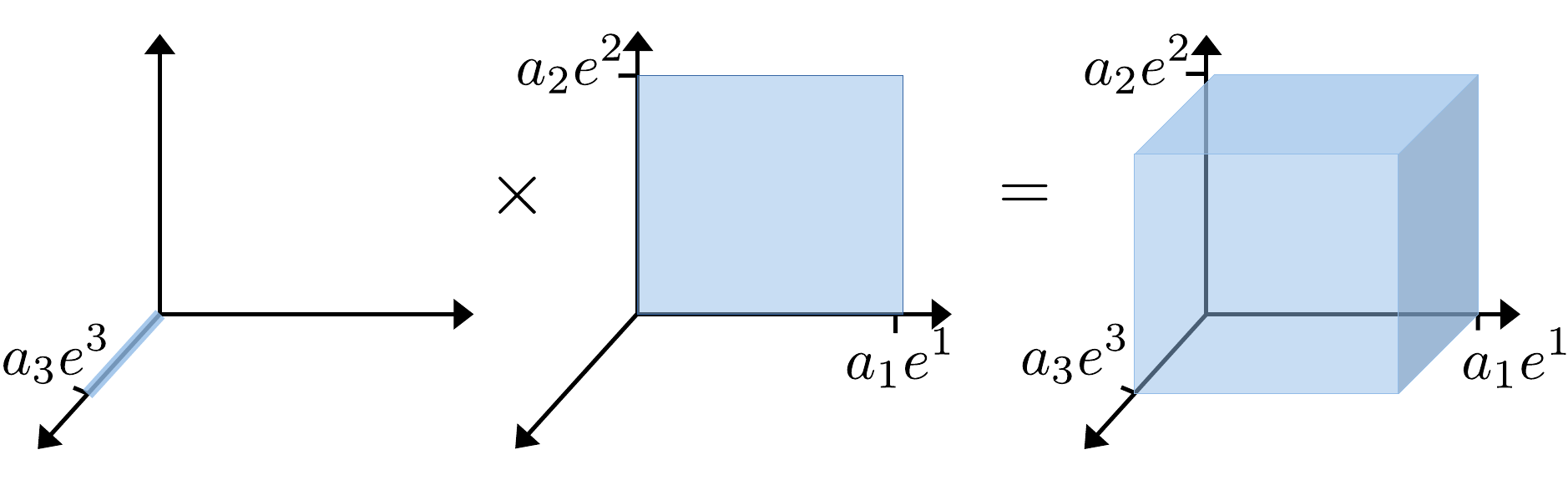}
  \caption{Multiplying
   the volume of a two-dimensional parallelotope and a one-dimensional parallelotope.}
    \label{fig:alpha_contr}
\end{center}
\end{figure}

\begin{Example}
Consider the case~$n=3$, $k=2$, and~$A=\diag(a_1,a_2,a_3)
$. Then
$
A^{(2)}=\diag(
  a_1 a_2 , a_1a_3 ,a_2a_3  )
$,
\be\label{eq:U3}
U=\begin{bmatrix}
   0&0&-1\\0&1&0\\-1&0&0
\end{bmatrix},
\ee
and~$
   U^{T}AU=\diag(a_3,a_2,a_1) $.
Thus,
\begin{align} \label{eq:a1a2a3}
 (A^{(2)})^TU^TA U= 
  (a_1a_2a_3)I_3,
\end{align}
and this agrees 
  with~\eqref{eq::A^(k)_A^(n-k)_alt}.
  
  To provide a geometric 
  intuition for~\eqref{eq:a1a2a3},
   let~$e^i$, $i=1,2,3$, denote the~$i$th canonical vector in~$\R^3$. Assume for simplicity that~$a_i\geq 0$, $i=1,2,3$. Then $(A^{(2)})_{11}$ is the volume of a parallelotope  
  with vertices~$a_1 e^1,a_2 e^2 $, and~$(U^TAU)_{11}$ is the volume of the parallelotope (in fact, line) with vertex~$a_3 e^3$. The product of these two volumes is the volume of a parallelotope with vertices~$a_1 e^1$, $a_2e^2$, and~$a_3e^3$, i.e., $\det(A)$
  (see Fig.~\ref{fig:alpha_contr}).
\end{Example}

 If~$n$ is even, then taking~$k=n/2$ in Thm.~\ref{thm:dual_,ult} yields the following result. 
 \begin{Corollary}
\label{coro:n/2}
Let~$A\in \C^{n\times n}$, with~$n$ even. Let~$r:=\binom{n}{n/2}$. Then
\begin{align}\label{eq:neven}
  ( U A^{(n/2)} )^T 
 A^{(n/2)}U=\det(A)I_r   . 
\end{align}
 \end{Corollary}

\subsection{Duality between additive compounds}
The next result describes a duality relation between the   additive  compound matrices~$A^{[k]}$ and~$A^{[n-k]}$.
\begin{Theorem}\label{thm:dual_add}
 Let~$A\in\C^{n\times n}$. Then
\begin{align}\label{eq:A[k]=I-UTA[n-k]U}
    (A^{[k]})^T+U ^TA^{[n-k]}U =\trace(A)I_r.
\end{align}
\end{Theorem}
\begin{IEEEproof}
Fix~$\varepsilon>0$. 
  Thm.~\ref{thm:dual_,ult}
yields
\begin{align}\label{eq:rpof}
    \big((I_n+\varepsilon A)^{(k)}\big)^TU^T(I_n+\varepsilon A)^{(n-k)}U
    &= \det(I+\varepsilon A)I_r.
\end{align}
By~\eqref{eq:(I+epsilonA)^k}, 
the term on the left-hand side   of~\eqref{eq:rpof} is
\begin{align*}
 (I_r+\varepsilon A^{[k]})^T & U^T(I_r+\varepsilon A^{[n-k]})U+o(\varepsilon)\\
    &=I_r+\varepsilon (A^{[k]})^T +\varepsilon U^TA^{[n-k]}U+o(\varepsilon).
\end{align*}
The term on the right-hand side of~\eqref{eq:rpof} is
\begin{align*}
\det(I+\varepsilon A)I_r=(1+\varepsilon \trace(A) +o(\varepsilon) )I_r.
\end{align*}
We conclude that
\[
 \varepsilon (A^{[k]})^T+\varepsilon U^TA^{[n-k]}U+o(\varepsilon)=(\varepsilon \trace(A) +o(\varepsilon) )I_r.
\]
Dividing both sides by~$\varepsilon$,  and taking~$\varepsilon \to 0$
completes the proof. 
\end{IEEEproof}

\begin{Example}
Consider the case~$n=3$, $k=2$, and~$A=\diag
(a_1, a_2,a_3)
$. Then
$
A^{[2]}=\diag(
  a_1+ a_2 ,
  a_1+a_3,
  a_2+a_3)  
$,  $U$ is as in~\eqref{eq:U3},
and~$
   U^{T}AU=
   \diag(
  a_3,
  a_2,
  a_1 )$.
Thus,
\begin{align*}
    (A^{[2]})^T+ 
    U ^TA U =(a_1+a_2+a_3) I_3.
\end{align*}
and this agrees 
  with~\eqref{eq:A[k]=I-UTA[n-k]U}. 
\end{Example}

\begin{Remark}
Ref.~\cite{Muldowney1998}
 includes a result that is similar to Thm.~\ref{thm:dual_add}. However, the result there uses a different matrix~$U$, and is in fact wrong. A counterexample to the result as  stated  in Ref.~\cite{Muldowney1998} is, for example, the case~$n=4$ and~$k=2$.  
\end{Remark}

\begin{Remark}
It follows from~\eqref{eq:A[k]=I-UTA[n-k]U}
that
$
  A^{[n-1]}
  =\trace(A)I_r -
   U ^TA^T  U   .
$
This special case  already appeared in~\cite{schwarz1970}, and has been used in the analysis of~$(n-1)$-positive systems~\cite{Eyal_k_posi}.
\end{Remark}

 One implication of
 Thm.~\ref{thm:dual_add}, that will be used below, is the following. 
 
  \begin{Corollary}\label{coro:commute}
The matrices~$(A^{[k]})^T$ and~$U^TA^{[n-k]}U$ commute.
\end{Corollary}
\begin{IEEEproof}
It follows from Thm.~\ref{thm:dual_add} that
$ (A^{[k]})^T=\trace(A)I_r-U^TA^{[n-k]}U $, so
\begin{align*}
    (A^{[k]})^T U^TA^{[n-k]}U&=(\trace(A)I_r-U^TA^{[n-k]}U)U^T A^{[n-k]}U\\
    &=U^TA^{[n-k]}U(\trace(A)I_r-U^TA^{[n-k]}U)\\
    &=U^TA^{[n-k]}U(A^{[k]})^T,
\end{align*}
and this completes the proof.
\end{IEEEproof}

  If~$n$ is even then taking~$k=n/2$ in Thm.~\ref{thm:dual_add} yields the following result. 
 \begin{Corollary}
\label{coro:add_n/2}
Let~$A\in \C^{n\times n}$, with~$n$ even. Let~$r:=\binom{n}{n/2}$. Then
\begin{align}\label{eq:add_neven}
     (A^{[n/2]})^T+U ^TA^{[n/2]}U =\trace(A)I_r.
\end{align}
 \end{Corollary}

 \subsection{Duality between multiplicative compounds of the matrix exponential}
 The next result uses the duality relations above to derive a duality relation for the multiplicative  compounds of the exponential of a matrix. 
\begin{Theorem}
 Let~$A\in\C^{n\times n}$. Then
 \begin{equation}\label{eq:e^(k)-e^(n-k)}
 ( (    \exp(A)   )^{(k)}  )^T=\exp(\trace(A) ) U^T\big(\exp(-A)\big)^{(n-k)}U.
\end{equation}
\end{Theorem}
\begin{IEEEproof}
Let~$D:= \big(  (\exp(A))     ^{(k)}\big)^T      U^T(\exp(A))^{(n-k)}U $. 
Using the identity
$    (\exp(A))^{(k)}=\exp(A^{[k]})
$
(see,  e.g.,~\cite{muldowney1990compound}) 
and  the fact that~$U^T=U^{-1}$
yields 
\begin{align*}
    U^T(\exp(A))^{(n-k)}U=
    U^T \exp(A^{[n-k]})U=\exp( U^TA^{[n-k]}U).
\end{align*}
Thus, $D=(\exp(  A^{[k]} ))^T 
    \exp( U^TA^{[n-k]}U) $,
and applying  Corollary~\ref{coro:commute}
gives
 \begin{align*}
  D&=   
     \exp( (A^{[k]})^T+U^TA^{[n-k]}U)\\
    &=\exp ( \trace(A)I_r ) \\
    &=\exp(\trace(A)) I_r,
\end{align*}
and this completes the proof.
\end{IEEEproof}
 
 \begin{Example}
 For~$n=3$ and~$k=2$,
 \eqref{eq:e^(k)-e^(n-k)} becomes
 \[
 ( (    \exp(A)   )^{(2)}  )^T=\exp(\trace(A))U^T 
 \exp(-A) U.
\]
This provides an expression for  all the~$2$-minors of~$\exp(A)$ that does not require to  compute any minors.
\end{Example}

\subsection{Duality between Log Norms of Additive Compounds}
\begin{Proposition}\label{prop:measure_inequality}
Let~$A\in\R^{n\times n}$. Fix~$k\in\{1,\hdots,n-1\}$, and let~$r=\binom{n}{k}$. Let~$U\in\{-1,0,1\}^{r\times r} $ 
be the matrix defined in~\eqref{eq:U_r}. Then for any log norm~$\mu:\R^{r\times r}\to\R$, we have:
\begin{equation}\label{eq:measure_inequality}
    \mu(A^{[k]}) =  \trace(A)+\mu_{U^T}(-(A^{[n-k]})^T).
\end{equation}
\end{Proposition}
\begin{IEEEproof}
By~\eqref{eq:A[k]=I-UTA[n-k]U},
\[
   A^{[k]}   =\trace(A)I_r-U ^T (A^{[n-k]})^T U.
\]
Applying~$\mu$ on both sides of this equation and using the fact that~$\mu(c I+B) =c+\mu (B) $ for any scalar~$c$  (see, e.g.,~\cite{vid_desoer})  yields:
\begin{align*}
    \mu(A^{[k]})&= \trace(A) +\mu ( -U^T (A^{[n-k]})^T U ),
\end{align*}
and this completes the proof. 
\end{IEEEproof}

\begin{Corollary}\label{coro:fulan/2}
Let~$A\in\R^{n\times n}$, with~$n$ even. Let~$r:=\binom{n}{n/2}$.
 Then for any log norm~$\mu:\R^{r\times r}\to\R$, we have:
\begin{equation}\label{eq:measure_inequality_n2}
    \mu(A^{[n/2]}) =  \trace(A)+\mu_{U^T}(-(A^{[n/2]})^T).
\end{equation}
In particular, for any log norm~$\mu_p$ induced by an~$L_p$ norm, we have
\begin{equation*} 
    \mu_p(A^{[n/2]}) =  \trace(A)+\mu_q(-(A^{[n/2]})),
\end{equation*}
and for the log norm  induced by the~$L_2$
norm:
\begin{equation}\label{eq:measure_inequality_n2_LP}
    \mu_2(A^{[n/2]}) =  \trace(A)+\mu_2(- A^{[n/2]} ). 
\end{equation}

\end{Corollary}

\begin{Example}
Consider~$A\in\R^{(2\ell)\times (2\ell)}$, with~$A =\diag(
\lambda_1,\dots,\lambda_{2\ell})$, and~$\lambda_1\geq\dots\geq\lambda_{2\ell}$. Then~$\tr(A)=\sum_{i=1}^{2\ell} \lambda_i$,
$ \mu_2(A^{[n/2]}) =\sum_{i=1}^{ \ell} \lambda_i $, and 
$ \mu_2(-A^{[n/2]}) = - \sum_{i= \ell+1}^{2 \ell} \lambda_i $, so clearly~\eqref{eq:measure_inequality_n2_LP} holds. 
\end{Example}

The formulas in  Prop.~\ref{prop:k_add_meas}
are not always easy to use. Indeed, the cardinality 
of the  set~$Q(k,n)$ is~$\binom{n}{k}$ and this may be very large.  
In the next section, we use the duality relation~\eqref{eq:measure_inequality}
to derive a sufficient condition for $k$-contraction of non-linear dynamical systems
which does not require calculating   compounds of the Jacobian.

\section{Compound-free sufficient condition for $k$-contraction}\label{sec:comp_free}
We begin by defining a new matrix operator.
\begin{Definition}\label{def:deftau}
Given an integer~$n\geq 1$, $ k\in\{1,\dots,n \}$,~$p\in\{1,2,\infty\}$, and an invertible matrix~$T\in\R^{n\times n}$,  the
\emph{$k$-shifted log norm}
$\tau_{p,k}:\R^{n\times n}\to\R$ is defined by
\be\label{eq:deftau}
\tau_{p,k}(A):=\tr(A) + (n-k) \mu_{q,T}(-A) , 
\ee
where~$q$ is such that~$q^{-1}+p^{-1}=1$.
\end{Definition}
Note that the name is justified by the equality $\tau_{p,k}(A) =\mu_{q,T} (\tr(A)I_n + (n-k) (-A))$. Note also that
the matrix~$U$ does not appear in the definition of~$\tau_{p,k}$.
As we will see below, this is 
  due to Lemma~\ref{lem:pnorm_U_inv},

We can now state the  main result in this section.  
 \begin{Theorem}\label{thm:k_cont_no_add}
 For any~$A\in\R^{n\times n}$, $ k\in\{1,\dots,n\}$,~$p\in\{1,2,\infty\}$, and an invertible matrix~$T\in\R^{n\times n}$, we have 
    \begin{equation}\label{eq:mu_k_add_ub}
       \mu_{p,T^{(k)}}(A^{[k]} ) \le  \tau_{p,k}(A).
    \end{equation}
\end{Theorem}
In other words, $\tau_{p,k}(A)$ provides an upper bound on~$\mu_{p,T^{(k)}}(A^{[k]} )$
that does not require to compute any compounds. 
    In particular, if
$   \tau_{p,k} (A)\le -\eta < 0$
    then~$\dot x=Ax$ is 
 $k$-contracting with rate~$\eta$ w.r.t.
 the scaled~$L_p$ norm with weight matrix~$T^{(k)}$.
Similarly, if the Jacobian~$J(t,x)$  of~\eqref{eq:sys}
satisfies
\[
\tau_{p,k} (J(t,x))\le -\eta < 0,\text{ for all }t\geq 0, x\in\Omega
\]
then~\eqref{eq:sys}
is 
 $k$-contracting with rate~$\eta$ w.r.t.
 the scaled~$L_p$ norm with weight matrix~$T^{(k)}$.
Note that for~$k=1$  this   
  yields  a non-standard  sufficient condition for contraction w.r.t.~$L_p$, namely, 
\[
\tr(J(t,x)) + (n-1) \mu_{q,T}(-J(t,x))
\leq-\eta<0
\]
for all~$t\geq 0$, $x\in\Omega$.  
More importantly, for~$k>1$, 
this provides a sufficient condition for~$k$-contraction that does not require computing any compounds.\footnote{More precisely, it requires to   compute only the trivial compounds~$J=J^{[1]}$ and~$\tr(J)=J^{[n]}$.}


The remainder of this section is devoted to the proof of Thm.~\ref{thm:k_cont_no_add}.  This requires the following
auxiliary result
that may be of independent interest.

\begin{Proposition}\label{prop:meas_ub}
    Fix~$A \in \R^{n \times n}$,
    $p \in \{1,2,\infty\}$, and~$k,\ell \in \{1,\dots,n\}$ with $\ell \le k$. Let~$T \in \R^{n \times n}$ be invertible. Then
    \begin{align}\label{eq:uku2}
         \frac{1}{k} \mu_{p,T^{(k)}}(A^{[k]}) &\le \frac{1}{\ell} \mu_{p,T^{(\ell)}}(A^{[\ell]}).
    \end{align}
\end{Proposition}

For example, for $\ell=1$ this gives   $\mu_{p,T^{(k)}}(A^{[k]}) \le k \mu_{p,T} (A)$.

\begin{IEEEproof}
We will use the following easy to verify fact. If~$a_1\geq\dots \geq a_n$ and~$k\in\{1,\dots,n-1\}$ then 
\be \label{eq:simptpver}
        \frac{1}{k}\sum_{i=1} ^k a_i - \frac{1}{k+1}\sum_{i=1} ^{k+1} a_i \geq \frac{  a_k-a_{k+1}}{k+1}\geq 0.
\ee

    We begin by proving~\eqref{eq:uku2} for~$T=I_n$.  
   We first consider the case~$p=2$. Let~$\lambda_1    \ge \dots \ge \lambda_n$ denote the eigenvalues of~$(A + A^T)/2$. Then 
    \begin{align*}
        \frac{\mu_2(A^{[k]}) }{k} -\frac{\mu_2(A^{[k+1]}) }{k+1}  &= \frac{1}{k}\sum_{i=1}^k \lambda_i  -\frac{1}{k+1}\sum_{i=1}^{k+1
        } \lambda_i  \\
&\geq 0 . 
 \end{align*}

   We now consider the case~$p=1$. 
   For any~$\alpha \in Q(k,n)$, 
   let $c_{\alpha,i} := a_{ii} + \sum_{j \notin \alpha} |a_{ji}|$, with~$i\in\{1,\dots,n\}$. 
   Let~$\beta := \argmax_{\alpha \in Q(k,n)} \sum_{i \in \alpha}c_{\alpha,i}$. Let $i_1,\dots,i_k \in \{1,\dots,n\}$ be such that $c_{\beta,i_1} \ge\dots \ge  c_{\beta,i_k}  $. Using Prop.~\ref{prop:k_add_meas} and~\eqref{eq:simptpver}
   gives 
    \begin{align*}
         {k}^{-1}\mu_1\left(A^{[k]}\right) &= \max_{\alpha \in Q(k,n)}  {k}^{-1}\sum_{i \in \alpha}  c_{\alpha,i}  \\
        &=  {k}^{-1} (c_{\beta,i_1} + \dots + c_{\beta,i_k}) \\
        &\le  {\ell}^{-1} (c_{\beta,i_1} + \dots + c_{\beta,i_\ell}) \\
        &\le \max_{\gamma \in Q(\ell,n)}  {\ell}^{-1} \max_{i \in \gamma}  ( a_{ii} + \sum_{\substack{j \notin \gamma}}  |a_{ji}|  ) \\
        &=  {\ell}^{-1} \mu_1(A^{[\ell]}).
    \end{align*}
    
    The proof for the  case~$p=\infty$ is similar and thus omitted.
    We conclude that
        \begin{align*} 
         \frac{1}{k} \mu_{p }(A^{[k]}) &\le \frac{1}{\ell} \mu_{p }(A^{[\ell]}).
    \end{align*}
To complete the proof of~\eqref{eq:uku2}, fix an  invertible matrix~$T$. Using~\eqref{eq:add_T} and the fact that~$\mu_{T}(A) = \mu(TAT^{-1})$, we have that
    \begin{align*}
      k^{-1} \mu_{p,T^{(k)}}(A^{[k]}) &=  k^{-1} \mu_p((TAT^{-1})^{[k]}) \\
        &\le  {\ell} ^{-1} \mu_p((TAT^{-1})^{[\ell]}) \\
        &= {\ell} ^{-1} \mu_{p,T^{(\ell)}}(A^{[\ell]}),
    \end{align*}
    and this completes the proof.
    \end{IEEEproof}


\begin{Example}
    Let~$A=I_n$ and fix~$k\in\{1,\dots,n\}$.
    Then~$A^{[k]}=k I_r$, so for any monotonic  norm~\footnote{Recall that a  norm~$|\cdot|:\R^n\to\R_{\ge0}$ is \emph{monotonic} if~$|y_i| \leq |x_i|$,   $i = 1,\dots,n$, implies 
    that~$|y| \leq |x|$. All $L_p$ norms are monotonic;  see~\cite{Bauer1961_mono_norms}.} 
    and any invertible matrix~$T$, we have~$\mu_{p,T^{(k)}}(A^{[k]}) = k$.
    Thus,   inequality~\eqref{eq:uku2} holds with an equality, implying that the bound  cannot be improved in general. 
\end{Example}

\begin{Remark}
Prop.~\ref{prop:meas_ub}
    implies in particular that if   the system~\eqref{eq:sys} 
    satisfies the infinitesimal condition for~$\ell$-contraction with rate~$\eta$ in~\eqref{eq:infik} w.r.t. an~$L_p$ norm scaled by~$T^{(\ell)}$, with~$p\in\{1,2,\infty\}$,
    then  for any~$k\in\{\ell,\ell+1,\dots,n\}$ it is~$k$-contracting 
    with rate~$\frac{k}{\ell}  \eta $ w.r.t. the same norm scaled by~$T^{(k)}$ 
     (see also~\cite{kordercont,Angeli2021nonosc}). More generally, if a system is~$\alpha$-contracting, with~$\alpha>0$ real, then 
    it is also~$(\alpha+\varepsilon)$-contracting  for any~$\varepsilon\geq 0$~\cite{wu2020generalization}.  
\end{Remark}

\begin{Remark}
Suppose that the system~\eqref{eq:sys}  satisfies the infinitesimal condition for~$k$-contraction w.r.t. to~$L_p$ for some~$p\in\{1,2,\infty\}$.
Then the system is also~$n$-contractive, that is,~$\tr(J(x))<0$ for all~$x\in\Omega$. 
Suppose now that either~$p=1$ or~$p=\infty$. 
Since  there exists at least one diagonal entry of~$J(x)$ that is negative,
\be\label{eq:munega}
\mu_p(-J(x))>0,\text{ for all } x\in\Omega.
\ee
For~$p=2$ we have
$\mu_2(-J(x))=\lambda_1( -(J(x)+J^T(x))/2 )$ and since
\[
\tr( -(J(x)+J^T(x))/2 )=\tr(-J(x))>0,
\]
the formula for the $L_2$ log norm implies that 
Eq.~\eqref{eq:munega} holds also when~$p=2$. 
Thus,  the sufficient condition for~$k$-contraction in Thm.~\ref{thm:k_cont_no_add} is a trade-off between the negativity of~$\tr(J(x))$ and the positivity of~$(n-k)\mu_p(-J(x))$.  
In this case, if the sufficient condition for $k$-contraction holds, that is,
\[
tr(J)+(n-k)\mu_p(-J)\leq-\eta<0 ,
\]
then clearly  
\[
tr(J)+(n-(k+1))\mu_p(-J)\leq-\eta<0 
\]
so the sufficient condition for~$(k+1)$-contraction also holds. 
\end{Remark}

 We can now prove Thm.~\ref{thm:k_cont_no_add}.
 
 \begin{IEEEproof}[Proof of Thm.~\ref{thm:k_cont_no_add}]
    Using~\eqref{eq:weight_mat_meas} and~\eqref{eq:add_T} yields 
    \begin{align*}
        \mu_{p,T^{(k)}}(A^{[k]}) &= \mu_p(T^{(k)}A^{[k]}(T^{(k)})^{-1}) \\
        &= \mu_p((TAT^{-1})^{[k]}) \\
        &= \mu_q(((TAT^{-1})^{[k]})^T),
    \end{align*}
    where we used the duality relation for log norms in Lemma~\ref{lem:mat_meas_dual}. Applying the duality   for additive compounds in Thm.~\ref{thm:dual_add}, and using the fact that~$\mu(A + cI) = \mu(A) + c$ and~$\mu(cI) = c$ for all log norms and any~$c \in \R$, we get
    \begin{align}\label{eq:softr}
        \mu_{p,T^{(k)}}(A^{[k]}) &= \tr(TAT^{-1}) + \mu_q(- U^T(TAT^{-1})^{[n-k]}U)\nonumber \\
        &= \tr(A) + \mu_q(-(TAT^{-1})^{[n-k]}) \nonumber\\
        &= \tr(A) + \mu_{q,T^{(n-k)}}(-A^{[n-k]}),
    \end{align}
    where we used Lemma~\ref{lem:pnorm_U_inv}. Applying Prop.~\ref{prop:meas_ub}   yields
    \begin{equation}
        \mu_{p,T^{(k)}}(A^{[k]}) \le \tr(A) + (n-k)\mu_{q,T}(-A),
    \end{equation}
    and this completes the proof.
\end{IEEEproof}

The next result summarizes some     properties of 
the $k$-shifted log norm~$\tau_{p,k}$.
The proof follows from Definition~\ref{def:deftau}, linearity of the trace operator, 
and known properties of log norms
(see, e.g.,~\cite{vid_desoer}).
\begin{Proposition}
Fix~$A,B\in\R^{n\times n}$, $k\in\{1,\dots,n\}$, and~$p\in\{1,2,\infty\}$. Then
\begin{enumerate}
    \item
    $\tau_{p,k}(0)=0 $.
\item
$|\tau_{p,k}(A ) -\tau_{p,k}(B )| \leq |\tr(A-B)|+(n-k) \|A-B\|_{q,T}$. 
\item $\tau_{p,k}(A+B) \leq  \tau_{p,k}(A)+\tau_{p,k}(B)$. 
\item $\tau_{p,k}(cA)=c\tau_{p,k}(
A)$, for any~$c\in\R_+$. 
\item $\tau_{p,k}(A+cI_n )= \tau_{p,k}(
A) +k c  $, for any~$c\in \R $. In particular,
\[
\tau_{p,k}(I_n)=k  ,  \quad \tau_{p,k}(-I_n)=-k .
\]
\item
\begin{align*}
    \tr(A)-(n-k)\|A \| _{q,T}
&\leq
\tr(A)-(n-k)\mu_{q,T}(A)
\\&\leq \tau_{p,k}(A)\\&\leq \tr(A)+(n-k)\|A \| _{q,T}.
\end{align*}
\item
$
\tau_{p,k}(rA+(1-r)B) \leq
r\tau_{p,k}(A) +
(1-r)\tau_{p,k}( B)$,  for any  $r\in[0,1]$. 
\end{enumerate}
\end{Proposition}

In particular,~$\tau_{p,k}$ is continuous, 
sub-additive, positively homogeneous of degree one, and convex. The latter property
implies that it is possible to verify the sufficient condition for~$k$-contraction in Thm.~\ref{thm:k_cont_no_add} for a polytope of dynamical systems by checking only the vertices of this polytope.


\section{Applications}
We now describe several applications of Thm.~\ref{thm:k_cont_no_add}. In particular, we show how it can be used to prove $k$-contraction in $n$-dimensional nonlinear systems without computing any~$k$-compounds. However, it is instructive to  begin with LTI systems. 

\subsection{$k$-contraction in LTI systems}
Consider the LTI system
\be\label{eq:lti}
\dot x(t)=A x(t),
\ee
with~$A\in\R^{n\times n}$.

Suppose first that~$A=\diag(\lambda_1,\dots,\lambda_n)$, with
\be\label{eq:order}
\lambda_1\geq\dots\geq \lambda_n.
\ee
Fix~$p\in\{1,2,\infty\}$. Then
\begin{align*}
   \tau_{p,k}(A)&=\tr(A)+(n-k)\mu_q(-A)\\
    &= \lambda_1+\dots+\lambda_n-(n-k)\lambda_n\\
    &=-(n-k-1) \lambda_n+\lambda_1+\dots+\lambda_{n-1}.
    \end{align*}
Combining this with~\eqref{eq:order}
implies that~$\tau_{p,k} (A)<0$ iff
\be\label{eq:connk}
\lambda_1+\dots+\lambda_{n-1}<(n-k-1)\lambda_n<0.
\ee
If~$k=n-1$ then this is equivalent to~$\lambda_1+\dots+\lambda_{n-1}<0 $ which is indeed a necessary and sufficient condition for~$(n-1)$-contraction of~\eqref{eq:lti}.
For~$k<n-1$ condition~\eqref{eq:connk}
requires that the sum of the first~$n-1$ eigenvalues is ``negative enough''
to guarantee~$k$-contraction. 


Now assume that~$A$ is not necessarily diagonal. 
  Let~$\lambda_1,\dots,\lambda_n$ denote the eigenvalues of~$A$, ordered such that~$\real (\lambda_1 )\ge \dots\ge \real (\lambda_n)$. 
  Then~\eqref{eq:deftau}
  with~$T=I_n$ yields
  \[
   \tau_{p,k}(A)=\tr(A) + (n-k) \mu_q (-A) .
   \]
Using the   bound $\mu(-A)\geq \real(- \lambda_n)$ gives 
  \begin{align*}
     \tau_{p,k} (A)&  
          \geq \sum_{i=1}^n \real (\lambda_i) + (n-k) \real (-\lambda_n) \\
         &\geq n\real (\lambda_n )-(n - k)\real (\lambda_n)\\
         &= k\real( \lambda_n),
  \end{align*}
   so in particular if~$\tau_{p,k}(A)<0$ then we must have that~$\real (\lambda_n) <0$. In other words,
   if  the sufficient $k$-contractivity condition   holds then~$A$ is not ``too unstable'' in the sense that
   it  has  at least one eigenvalue with a  negative real part. 

Recall that for the LTI  system~\eqref{eq:lti} $k$-contraction implies that every   sum of~$k$ eigenvalues  of~$A$ has a negative real part~\cite{kordercont}.
Combining this with Thm.~\ref{thm:k_cont_no_add} yields the following result.
\begin{Corollary}
 Let~$A\in\R^{n\times n}$. Suppose that there
 exist an invertible  matrix~$T \in \R^{n \times n}$, 
    $k \in \{1,\dots,n\}$, and~$ q \in \{1,2,\infty\}$   such that  
\[
\tr(A) + (n-k) \mu_{q,T}(-A)< 0 . 
\]
Then every sum of~$k$ eigenvalues  of~$A$ has a negative real part.
\end{Corollary}

 The next result shows how Thm.~\ref{thm:k_cont_no_add} can be used to derive a  ``$k$-trace dominance condition'' guaranteeing $k$-contraction.
 
 \begin{Corollary}\label{eq:trac_domi}
Fix~$k\in\{1,\dots,n\}$. If there exist~$d_1,\dots,d_n>0$ such that
 \be\label{eq:scdom}
 -( n-k-1 )a_{ii}+\sum_{j\not =i}  \Big(a_{jj}+(n-k) \frac{d_j}{d_i} |a_{ji}  |\Big) \leq -\eta < 0,
 \ee
   for all~$i\in\{1,\dots,n\}$, 
 then
the LTI~\eqref{eq:lti} is~$k$-contractive with rate~$\eta$
w.r.t. the scaled  norm~$|\cdot|_{\infty,D}$, where~$D:=\diag(d_1,\dots,d_n)$.
 \end{Corollary}
 Note that for~$k=n$  condition~\eqref{eq:scdom} becomes
 \[
 \tr(A) \leq-\eta<0,
 \]
 whereas for~$k=n-1$ and~$D=I_n$
 it becomes
 \[
\tr(A)-a_{pp} + \sum_{j\not =p}   |a_{jp}  |
\leq-\eta<0, \text{ for all } p\in\{1,\dots,n\} .
 \]

\begin{IEEEproof}
We   prove Corollary~\ref{eq:trac_domi}
for the case~$D=I_n$. The general case follows by  replacing~$A$ with~$DAD^{-1}$. Consider 
\begin{align*}
\tau_{\infty,k} (A)&=\tr(A)+(n-k)\mu_1(-A)\\
&=\sum_{i=1}^n a_{ii}+(n-k) \max(c_1,\dots,c_n),
\end{align*}
where
$
c_i :=-a_{ii} +  
\sum_{j\not =i} |a_{ji}  |   
$, that is,
the sum of the entries in column~$i$ of~$(-A)$, with  off-diagonal entries taken with absolute value. 
For concreteness, assume that~$\max(c_1,\dots,c_n)=c_1$, that is,
\[
-a_{11} +  
\sum_{j\not =1} |a_{j1}  |\geq -a_{\ell\ell} +  
\sum_{j\not =\ell} |a_{j\ell}  |, \text{ for all } \ell \geq 1. 
\]
Then
\begin{align*}
\tau_{\infty,k} (A)& =\sum_{i=1}^n a_{ii}+(n-k) (-a_{11} +  
\sum_{j\not =1} |a_{j1}  |) \\
&= 
-( n-k-1 )a_{11}+\sum_{j\not =1}  (a_{jj}+(n-k)  |a_{j1}  |).
\end{align*}
Comparing this with~\eqref{eq:scdom}
completes the proof. 
\end{IEEEproof}

\subsection{$k$-contraction in LTV systems}

Applying  Thm.~\ref{thm:k_cont_no_add}  to prove~$k$-contraction 
requires    bounding~$\mu_{q,T}(-J(t,x))$ from above.
For the case of an LTV system and the~$L_2$ norm (i.e.,~$q=2$), a useful bound was derived in~\cite{Smith1986haus}.
\begin{Lemma}\cite{Smith1986haus}\label{lem:thetabo}
 Consider the matrix LTV  system
    \begin{equation}\label{eq:ltv}
        \dot X(t) = A(t) X(t), \quad X(0) = I,
    \end{equation}
    with~$A  : \R_{+ } \to \R^{n \times n}$   continuous. 
    Suppose that there exist~$Q  \succ 0$ and a continuous function~$\theta : \R_{+} \to \R$ such that
    \begin{equation}\label{eq:QA}
        A^T(t)Q + QA(t) + 2\theta(t) Q \succeq 0, \text{ for all } t\geq 0. 
    \end{equation}
Let~$P\succ 0$ be such that~$P^2=Q$. Then
\be\label{eq:mu2pbo}
   \mu_{2, P }(-A(t)) \le \theta(t),\text{ for all } t\geq 0  .
   \ee
\end{Lemma}
The proof follows from multiplying~\eqref{eq:QA} by~$P^{-1}$ on the left- and right-hand sides. 

Combining 
this bound with Thm.~\ref{thm:k_cont_no_add}
yields the following result.
\begin{Proposition}\label{prop:ltv_suff}
 Consider the matrix LTV  system~\eqref{eq:ltv}
  and suppose that the conditions in   Lemma~\ref{lem:thetabo} hold. 
  Fix~$k\in\{1,\dots,n\}$. If
\[
\tr(A(t)) +(n-k)\theta(t)\leq-\eta<0,  \text{ for all } t\geq 0 
\]
then the LTV system is~$k$-contracting with rate~$\eta$ w.r.t.   the scaled~$L_2$ norm~$|\cdot|_{2,P}$.
\end{Proposition}

\begin{IEEEproof}
Consider 
 \[
\tau_{2,k}(A(t))=\tr(A(t)) + (n-k) \mu_{2,P}(-A(t)) . 
 \]   
Combining this with~\eqref{eq:mu2pbo} gives
 \[
\tau_{2,k}(A(t))\leq \tr(A(t)) + (n-k) \theta(t) , 
 \]   
and applying Thm.~\ref{thm:k_cont_no_add} completes the proof. 
\end{IEEEproof}

\begin{Example}
Consider the LTI~\eqref{eq:lti}, but  with an uncertainty in the matrix~$A$.
A standard approach for modeling this is to assume that~$A$ is constant, unknown, and belongs to the convex hull of a set of~$s$ known matrices~$A_1,\dots,A_s$. We assume that all the matrices have the same trace
\[
r:=\tr(A_i),\quad  i=1,\dots,s.
\]
This is a typical case, for example, in modeling  biological interaction networks  (also known as chemical reaction networks), see, e.g.,~\cite{Angeli2021nonosc}.
Prop.~\ref{prop:ltv_suff}
implies that if we can find~$Q  \succ 0$ and~$\theta \in  \R$ such that
\[
A_i^TQ + QA_i + 2\theta Q \succeq 0, \text{ for all } i\in\{1,\dots,s\} ,
\]
and
\[
r +(n-k)\theta \leq-\eta<0,
\]
then the uncertain  LTI is $k$-contractive. We emphasize  again that this does not require computing any compounds. 
\end{Example}

\subsection{$k$-contraction in $n$-dimensional   Hopfield Neural Networks}
Consider the Hopfield neural network~\cite{hopfield_net}
\be\label{eq:hopf}
\dot x_i(t)=-\frac{x_i(t)}{r_i}+\sum_{j=1}^n w_{ij} \phi_j(x_j(t))+u_i    ,\quad i=1,\dots,n,
\ee
where~$r_i>0$, $u_i$ is a constant input to   neuron~$i$, $\phi_j:\R\to\R$ is the activation  function of neuron~$j$, and~$W=\{w_{ij}\}_{i,j=1}^n$ is the network
connection matrix. We assume that every~$\phi_i$ is~$C^1$. 

The stability of~\eqref{eq:hopf} has been studied extensively, e.g., via
Lyapunov analysis in \cite{ANM-JAF-WP:89,EK-AB:94,MF-AT:95}.
Ref.~\cite{YF-TGK:96} seems to be  the  first application of log norms to analyze Hopfield neural networks; later works include~\cite{cont_hopfield,AD-AVP-FB:21k} on contractivity w.r.t.  non-Euclidean
norms, and~\cite{MR-RW-IRM:20,LK-ME-JJES:21} on contractivity
w.r.t. Euclidean norms.
However, in many applications the network admits more than a single
equilibrium. For example, in using a Hopfield network as an associative
memory~\cite{hopfield_net,krotov2016}, every stored pattern  corresponds to an
equilibrium. Thus, if the network stores more than a single memory than it
cannot be contractive w.r.t. any norm.

The Jacobian  of~\eqref{eq:hopf} is  
\be\label{eq:jaco_hop}
J (x)=-\diag\big(r_1^{-1},\dots,r_n^{-1})+W\diag(\phi_1'(x_1) ,\dots,\phi_n' (x_n)\big),
\ee
where~$\phi_i'(x_i):=\frac{d}{dx}\phi_i(x)|_{x=x_i}$. 
Thm.~\ref{thm:k_cont_no_add} allows to derive a sufficient condition for~$k$-contraction of the Hopfield
network without computing compounds. One possibility is to assume that the~$\phi_i'$s are bounded,   and then apply  
the same approach as in 
Corollary~\ref{eq:trac_domi}.

\begin{Proposition}\label{prop:hopf}
Consider the Hopfield  network~\eqref{eq:hopf}.
Assume that the neuron activation
functions  satisfy
\be\label{eq:fibou}
0\leq m_i\leq|\phi_i'(z)|\leq M_i, \text{ for all } z\in\R \text{ and  } i\in\{1,\dots,n\}.
\ee
If there exist~$d_1,\dots,d_n>0$
such that
 \begin{align}\label{eq:scdom_hop}
 &-( n-k-1 ) (-r_i^{-1}  -  m_i|w_{ii}|)\nonumber \\
 &+\sum_{j\not =i}  
 \Big(-r_j^{-1} +M_j|w_{jj}| +(n-k) \frac{d_j}{d_i} M_i|w_{ji}| \Big) \leq-\eta<0,  
 \end{align}
  for all~$i\in\{1,\dots,n\}$,
 then~\eqref{eq:hopf} is $k$-contractive with rate~$\eta$ w.r.t. the scaled norm~$|\cdot|_{\infty,D}$, with~$D:=\diag(d_1,\dots,d_n)$.
\end{Proposition}
  
A common choice for the activation functions
in neural network models is~$\phi_i(z)=a_i\tanh(b_i z)$ and then clearly   condition~\eqref{eq:fibou}  indeed holds. Note also that if we set~$r_1=\dots=r_n=r$  then~\eqref{eq:scdom_hop} will hold 
for any~$r>0$ sufficiently small.
This makes sense, as
a smaller~$r$ makes~\eqref{eq:hopf}
``more stable''. We emphasize again that condition~\eqref{eq:scdom_hop} does not require to compute any compounds of the Jacobian~$J(x)$ in~\eqref{eq:jaco_hop}. 

\begin{IEEEproof}
Consider
\begin{align*}
\tau_{\infty,k}( J(x) ) (J(x))&=\tr(J(x))+(n-k)\mu_1(-J(x))\\
&=\sum_{i=1}^n(  w_{ii}\phi_i'(x_i)   -  r_i^{-1} )\\&+(n-k) \max (c_1(x),\dots,c_n(x)),
\end{align*}
where
\[
c_i(x):=r_i^{-1}-w_{ii}\phi'_i + \sum_{j\not =i} |w_{ji} \phi'_i | . 
\]
For concreteness, assume that~$\max(c_1(x),\dots,c_n(x))=c_1(x)$. Then
\begin{align*}
\tau_{\infty,k}( J(x) ) &= -(n-k-1) (w_{11}\phi_1'(x_1)   -  r_1^{-1} ) 
\\&+ \sum_{i\not =1} (  w_{ii}\phi_i'(x_i)   -  r_i^{-1} )
+(n-k)\sum_{j\not =1} |w_{ji} \phi'_i |.
\end{align*}
Applying~\eqref{eq:fibou} gives
\begin{align*}
\tau_{\infty,k}( J(x) )& \leq-(n-k-1) (-|w_{11} | m_1  -  r_1^{-1} ) 
\\
&+ \sum_{i\not =1} ( | w_{ii}|M_i   -  r_i^{-1} ) +(n-k)\sum_{j\not =1} |w_{ji}   |M_i
\end{align*}
for all~$x$, and this completes the proof. 
\end{IEEEproof}

\begin{figure}
 \begin{center}
  \includegraphics[scale=0.5]{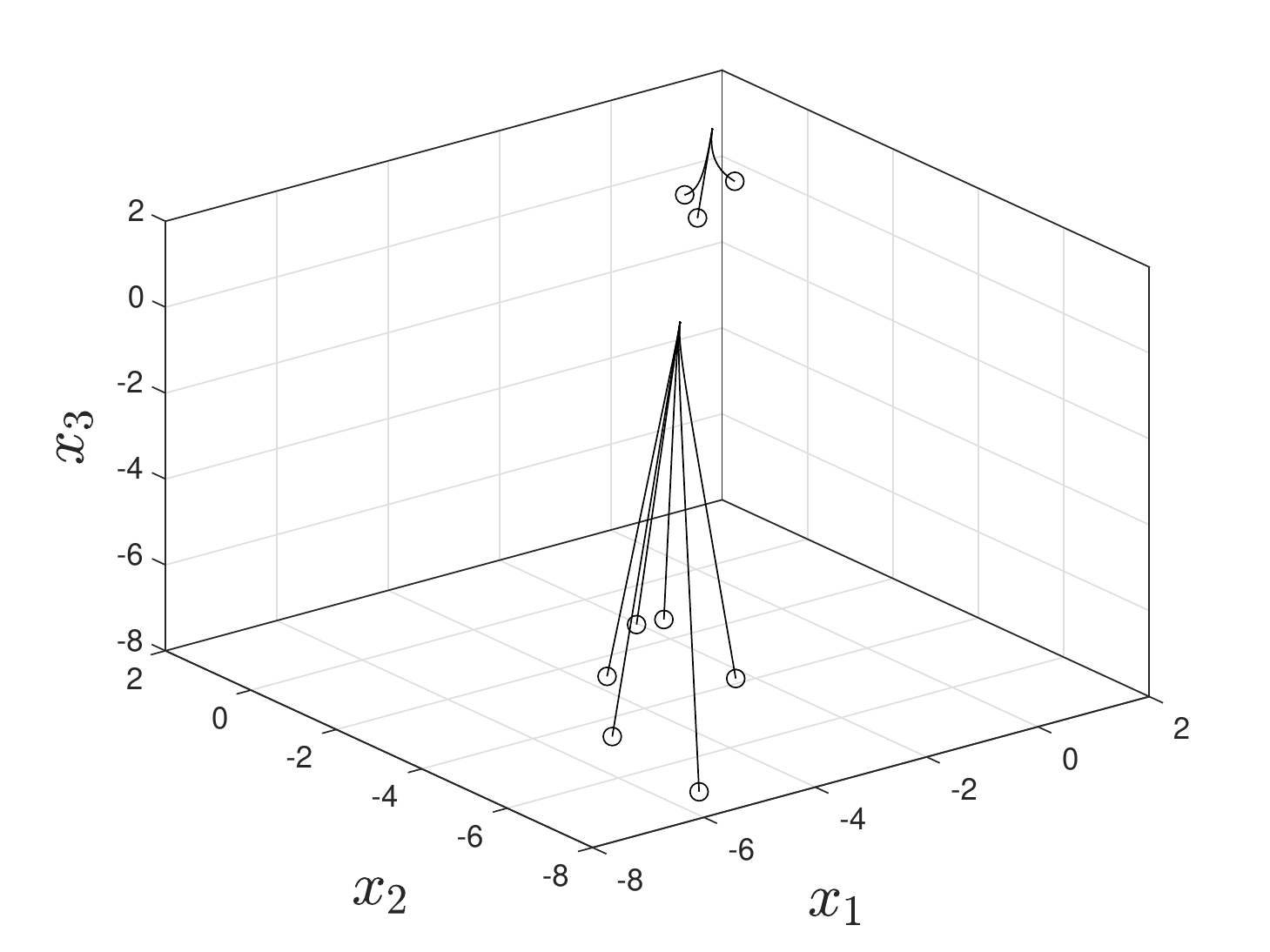}
	\caption{Trajectories of the Hopfield network in~\eqref{eq:hop3}.
	Initial conditions are marked by~o.}
	\label{fig:hopf} 
\end{center}
\end{figure}

\begin{Example}
Consider~\eqref{eq:hopf}
with~$|w_{ij}| = 1$ for all~$i,j$ (i.e., a binary weight matrix), $r_i=r$ and~$\phi_i(z)=\tanh(z)$ for all~$i$, that is, 
\be\label{eq:4Dhopf}
\dot x_i =-\frac{x_i}{r }+\sum_{j=1}^n (\pm 1)  \tanh(x_j )+u_i    ,\quad i=1,\dots,n ,
\ee
where~$\pm 1$ indicates a value that can be either~$-1$ or~$+1$. 
For large $r$, this system may have more than a single equilibrium and may not be contractive w.r.t. any norm. 
We apply Prop.~\ref{prop:hopf} with~$D=I_n$
to find a sufficient condition for $k$-contraction.
In this case,~$m_i=0$ and~$M_i=1$ for all~$i$, and 
\eqref{eq:scdom_hop} becomes 
 \[ 
  -k r ^{-1}    + (n-1)(n-k+1)\leq-\eta<0.
 \]
Thus, a sufficient condition for $k$-contraction is 
\be\label{eq:suff_kcon_hop}
r<\frac{k}{(n-1)(n-k+1)}.
\ee
Note also that this did not require to compute and analyze~$J^{[k]}(x) $ which in this case is an~$\binom{n}{k}\times \binom{n}{k}$
state-dependent
  matrix. 

For example,
if we require~$(n-1)$-contraction then~\eqref{eq:suff_kcon_hop}  becomes the condition~$r<1/2$.
 As a specific example, take~$n=3$, $w_{ij}=1$ for all~$i,j$, $r_i=0.49$ 
 and~$f_i(z)=\tanh(z)$ for all~$i$, that is,
 \be\label{eq:hop3}
 \dot x_i=-\frac{x_1}{0.49} +\sum_{i=1}^3 \tanh(x_i),\quad i=1,2,3.
 \ee
 This system admits at least three equilibrium points:
 \be\label{eq:jeq}
 e^1=0 , \; 
 e^2\approx 1.2447 \begin{bmatrix}
 1&1&1
 \end{bmatrix}^T 
 ,\;
 e^3=-e^2,
 \ee
 so it is not contractive w.r.t any norm. It satisfies the
  sufficient condition for~$2$-contraction (namely,~$r<1/2$), so the results of 
  Muldowney and Li~\cite{muldowney1990compound,li1995} imply that 
  every bounded solution converges to an equilibrium point. It is clear from~\eqref{eq:hop3} that all trajectories of the system are bounded. Fig.~\ref{fig:hopf} 
  depicts several trajectories of this system from random initial conditions. It may be seen that all the   trajectories indeed   converge to either~$e^2$ or~$e^3$.  
  Note that in this example, $n=3$ and~$k=2$, so we can easily compute and analyze~$J^{[2]}(x)$ directly, but our goal is merely 
  to demonstrate the bound  derived in Thm.~\ref{thm:k_cont_no_add}.
 \end{Example}

\subsection{Local stability of an equilibrium of a non-linear dynamical system}
Li and Wang~\cite{LI_WANG_STAB} proved that   a matrix~$A\in\R^{n\times n}$ is Hurwitz  iff the following two conditions hold:
\be\label{eq:twocond}
A^{[2]} \text{ is Hurwitz and } (-1)^n\det(A)>0.
\ee
 The proof is based on the fact that every eigenvalue of~$A^{[2]}$ is the sum of two eigenvalues  of~$A$. 
 This result was applied to prove the local stability of the endemic equilibrium~$e$ in an SEIR model with a  varying total population size~\cite{LI_WANG_STAB} by verifying that~$J(e)$, the Jacobian of the vector field evaluated at the equilibrium, is Hurwitz. In this case,~$J(e)\in\R^{3\times 3}$
 depends on various 
  parameters of the model and verifying  that~$J(e)$ is Hurwitz using the Routh–Hurwitz stability criterion is non-trivial. However, in general 
  verifying that~\eqref{eq:twocond} holds requires computing~$A^{[2]}$ which is a matrix of dimensions~$\binom{n}{2}\times\binom{n}{2}$. The next result uses the operator~$\tau_{p,k}$
  and does not require to compute~$2$-compounds.

 \begin{Corollary}
 Let~$e\in\R^n$ be an equilibrium of the system~$\dot x =f(x)$, with~$f\in C^1$. Let~$J(x):=\frac{\partial}{\partial x}f(x)$. 
 If there exist~$p\in\{1,2,\infty\}$ and an invertible matrix~$T\in\R^{n\times n} $ such that
 \be\label{eq:jep2}
 \tau_{p,2}(J(e))=\tr(J(e))+(n-2)\mu_{q,T}(-J(e)) <0,
 \ee
 where~$p^{-1}+q^{-1}=1$,
 and
 \be\label{eq:m1n}
 (-1)^n\det(J(e))>0
 \ee
 then~$e$ is locally asymptotically stable. 
 \end{Corollary}
\begin{IEEEproof}
By Thm.~\ref{thm:k_cont_no_add}, Eq.~\eqref{eq:jep2} implies that~$\mu( (J(e) )^{[2]} ) <0$, and combining this with~\eqref{eq:m1n} implies that~$J(e)$ is Hurwitz.
\end{IEEEproof}

Note that conditions~\eqref{eq:jep2} and~\eqref{eq:m1n} do not require to compute~$(J(e) )^{[2]}$. 

As a simple example consider again the Hopfield network in~\eqref{eq:hop3}
and the equilibrium points in~\eqref{eq:jeq}. We already know that condition~\eqref{eq:jep2} holds at any point in the state-space, so we only need to check condition~\eqref{eq:m1n}, that is,
\be  \label{eq:detcju}
\det(J(e))<0. 
\ee
 Using~\eqref{eq:jaco_hop} gives
 \begin{align*}
                J(x)&=-\diag(1/0.49,1/0.49,1/0.49)\\
                &+\begin{bmatrix}
                1-\tanh^2(x_1) &   1-\tanh^2(x_2) &  1-\tanh^2(x_3)\\
                1-\tanh^2(x_1) &   1-\tanh^2(x_2) &  1-\tanh^2(x_3)\\
                1-\tanh^2(x_1) &   1-\tanh^2(x_2) &  1-\tanh^2(x_3) 
                \end{bmatrix}.
 \end{align*}
 It follows that~$\det(J(x))=(-1/0.49)^2 ((-1/0.49) + 3-\sum_{i=1}^3  
 \tanh^2(x_i)
 ) $. In particular,~$\sign(\det(J(e^1)))>0$
 and~$ \sign(\det(J(e^2)))<0$. 
  Thus,~$e^2$ is locally  asymptotically stable,
  and~$e^1$ is not locally  asymptotically stable. 
 
\section{Conclusion}
Contraction theory plays an important role in
systems  and control theory. 
However, many systems cannot be analyzed using contraction theory. For example, systems with more than a single equilibrium are not contractive w.r.t. any norm.

The notion of $k$-contraction  provides a useful 
geometric generalization  of contraction theory, but
the standard sufficient 
condition for $k$-contraction of  $n$-dimensional systems 
may be difficult to verify, as it is based on a compound  matrix with dimensions~$\binom{n}{k}\times\binom{n}{k}$.
We derived    duality relations  between compound matrices, and used these to
develop a sufficient condition for  $k$-contraction
that does not require to  
 compute any compounds.  
We demonstrated  our approach  
by deriving a sufficient condition  for 
$k$-contraction of a Hopfield neural network.
In the particular case where~$k=2$ this implies that every bounded solution of the network converges to an equilibrium, which is of course a useful property when using the network as an associative memory~\cite{krotov2016}. 
We believe that the sufficient conditions for~$k$-contraction derived here will prove useful in
more   applications.

Several results in this paper
are proved  for~$L_p$ norms, with~$p\in\{1,2,\infty\}$. 
It may be of interest to 
try and generalize the proofs to any~$L_p$ norm. 
Another interesting direction for future research is to extend the tools described here to control synthesis. In other words, to systematically design a controller such that 
the closed-loop system  satisfies the  sufficient condition for~$k$-contraction.

 \begin{IEEEbiography}[{\includegraphics[width=1.1in,height=2in,clip,keepaspectratio]{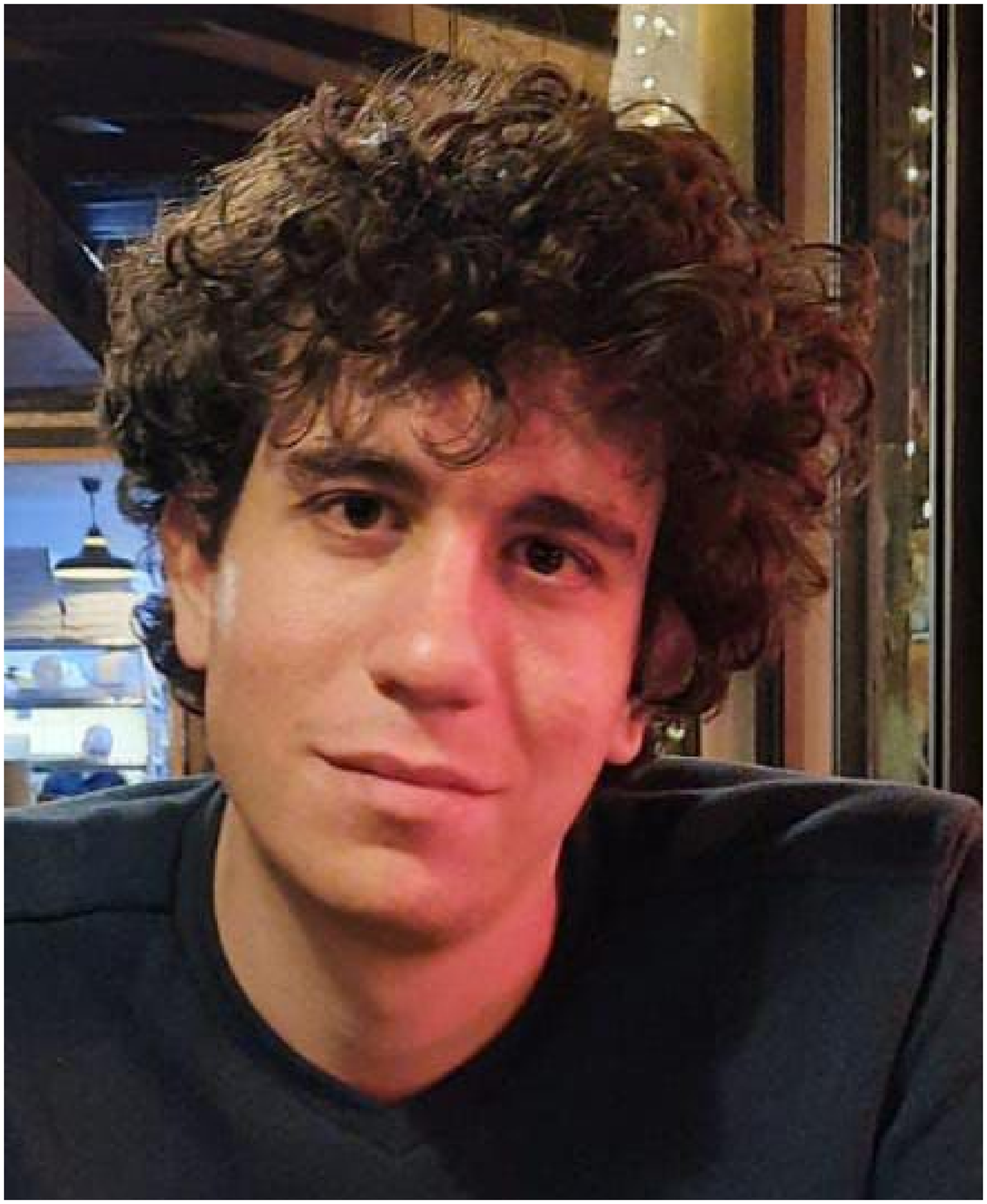}}\\]{Omri Dalin} received the  BSc degree in Mech. Eng. from Tel Aviv University, in 2020.
He is currently an MSc student at the Dept. of Mech. Eng. at Tel Aviv University. 
\end{IEEEbiography}

  \begin{IEEEbiography}
 {Ron Ofir} (Student Member, IEEE) received his BSc degree (cum laude) in Elec.  Eng. from the Technion-Israel Institute of Technology in 2019. He is currently pursuing his Ph.D. degree in the Vitebri Faculty of Electrical Engineering, Technion- Israel Institute of Technology. His current research interests include contraction theory and its applications in dynamics and control of power systems.
 \end{IEEEbiography}

\begin{IEEEbiography}[{\includegraphics[width=1in,height=1.5in,clip,keepaspectratio]{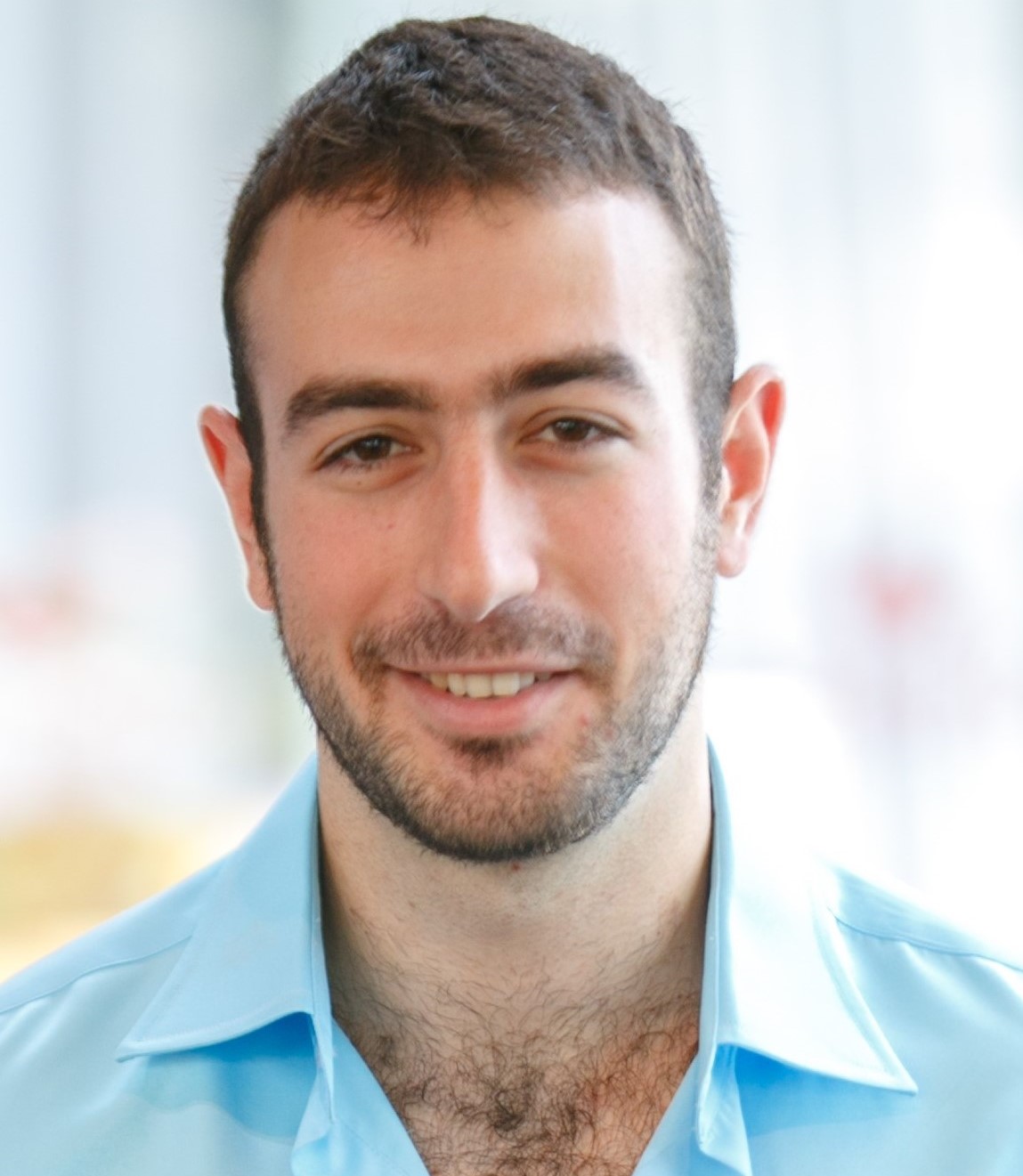}}]{Eyal Bar-Shalom} received the BSc degree (cum laude) and MSc degree (cum laude) in Elec. Eng. from Tel Aviv University, in 2012 and 2019, respectively.
He is currently a PhD student at the Dept. of Elec. Eng. at Tel Aviv University. 
\end{IEEEbiography}

 \begin{IEEEbiography}
    [{\includegraphics[width=1in,height=1.25in,clip,keepaspectratio]{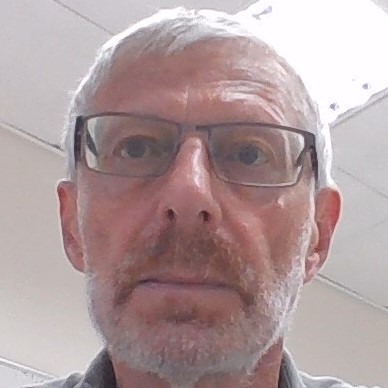}}]{ Alexander Ovseevich}
 received his MSc in Mathematics from the Moscow State University in 1973, and PhD in Mathematics from the Leningrad division of the Steklov mathematical institute, RAS in 1976. He got a degree of the doctor of Phys.-Math. Sciences from the Institute for Problems in Mechanics, Moscow, Russia in 1997. He worked in the Institute for Problems in Mechanics from 1978 to 2022 as a senior research fellow and a leading research fellow. In 2022 he joined the Dept. of Elec. Eng.-Systems, Tel Aviv University, where he is currently a research assistant. His research interests include many topics from Number Theory to Mathematical Physics and Control Theory. 
 \vspace*{-1.2cm}
\end{IEEEbiography}

\begin{IEEEbiography}
	[{\includegraphics[width=1in,height=1.25in,clip,keepaspectratio]{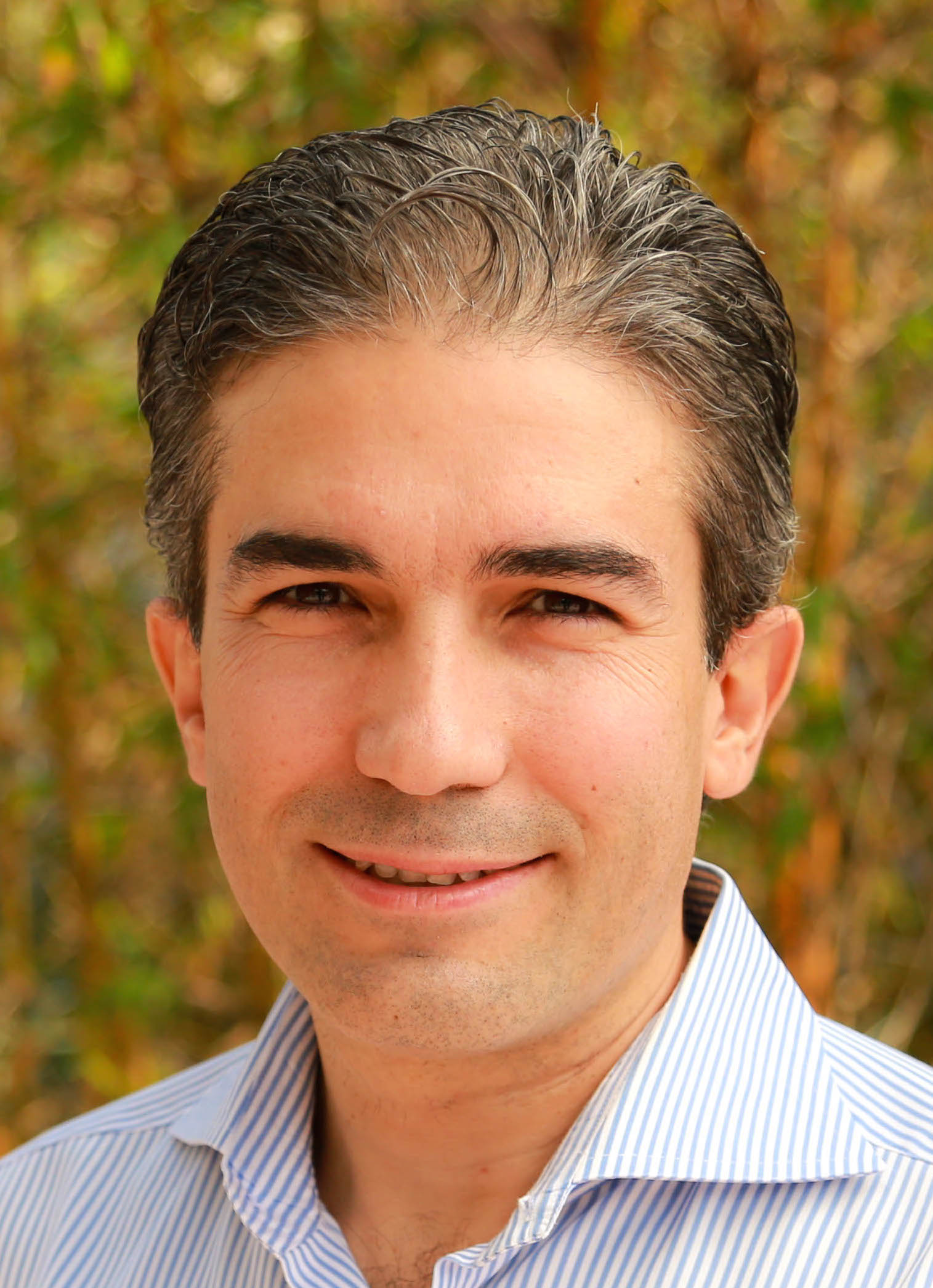}}]{Francesco Bullo}
	(Fellow, IEEE) is a
        Distinguished Professor of Mechanical Engineering at the University
        of California, Santa Barbara. He was previously associated with the
        University of Padova, the California Institute of Technology, and
        the University of Illinois. He served as IEEE CSS President and as
        SIAG CST Chair. His research focuses on contraction theory, network
        systems and distributed control with application to machine
        learning, power grids, social networks, and robotics. He is the
        coauthor of “Geometric Control of Mechanical Systems” (Springer,
        2004), “Distributed Control of Robotic Networks” (Princeton, 2009),
        “Lectures on Network Systems” (KDP, 2022, v1.6), and "Contraction
        Theory for Dynamical Systems" (KDP, 2022, v1.0).  He is a Fellow of
        ASME, IFAC, and SIAM.
\end{IEEEbiography}
 
\begin{IEEEbiography}
    [{\includegraphics[width=1in,height=1.25in,clip,keepaspectratio]{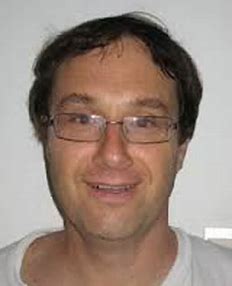}}]{Michael Margaliot}
 received the BSc (cum laude) and MSc degrees in
 Elec. Eng. from the Technion-Israel Institute of Technology-in
 1992 and 1995, respectively, and the PhD degree (summa cum laude) from Tel
 Aviv University in 1999. He was a post-doctoral fellow in the Dept. of
 Theoretical Math. at the Weizmann Institute of Science. In 2000, he
 joined the Dept. of Elec. Eng.-Systems, Tel Aviv University,
 where he is currently a Professor. His  research
 interests include the stability analysis of differential inclusions and
 switched systems, optimal control theory, computation with
 words, Boolean control networks, contraction theory, $k$-positive systems, and systems biology.
 He is co-author of \emph{New Approaches to Fuzzy Modeling and Control: Design and
 Analysis}, World Scientific,~2000 and of \emph{Knowledge-Based Neurocomputing}, Springer,~2009. 
 He  served
as  an Associate Editor of~\emph{IEEE Transactions on Automatic Control} during 2015-2017.
\vspace*{-1.2cm}
\end{IEEEbiography}

\end{document}